# ASYMPTOTICS OF SURVIVAL PROBABILITIES AND LOWER TAIL PROBABILITY PROBLEM

SVETLANA BOYARCHENKO AND SERGEI LEVENDORSKIĬ

ABSTRACT. The present paper is an addendum to the paper "Lévy models amenable to efficient calculations", where we introduced a general class of Stieltjes-Lévy processes (SL-processes) and signed SL processes defined in terms of certain Stieltjes-Lévy measures. We demonstrated that SL-processes enjoyed all properties that we used earlier to develop efficient methods for evaluation of expectations of functions of a Lévy process and its extremum processes, and proved that essentially all popular classes of Lévy processes are SL-processes; sSL-processes fail to possess one important property. In the present paper, we use the properties of (s)SL-processes to derive new formulas for the Wiener-Hopf factors $\phi_q^\pm$ for small $q$ in terms of the absolute continuous components of SL-measures and their densities, and calculate the leading terms of the survival probability also in terms of the absolute continuous components of SL-measures and their densities. The lower tail probability is calculated for more general classes of SINH-regular processes constructed earlier.

KEY WORDS: Survival probability, lower tail probability problem, Wiener-Hopf factorization, Stieltjes-Lévy processes, SINH-regular Lévy processes, KoBoL, Normal inverse Gaussian processes, Variance Gamma processes, stable Lévy processes

MSC2010 CODES: 60G51, 60G52,60-08,65C05,91G05,91G20,97M30

## 1. INTRODUCTION

Let $X$ be a Lévy process on $\mathbb{R}$ and $\bar{X}, \underline{X}$ the supremum and infimum processes of $X$, all starting at 0. In the paper, we derive asymptotic formulas for the survival probability $\mathbb{P}[\bar{X}_t < x]$, where $x > 0$ is fixed and $t \to \infty$, and lower tail probability, when $t$ is fixed and $x \downarrow 0$. In the case of stable Lévy processes, the problems are equivalent but this is not the case for other classes of processes. Both problems appear in various branches of natural sciences, finance and insurance, and have been extensively studied (see, e.g., [26, 3] and the bibliographies therein). However, the extant results do not cover all classes of Lévy processes. We consider wide classes of Lévy processes with exponentially decaying Lévy densities, and outline extensions of the results to the lower tail probability problem for mixtures of stable Lévy processes and processes with exponentially decaying Lévy densities; extension of the results for the survival function to this case remains an open problem.

The proofs are based on the integral representations for $\mathbb{P}[\bar{X}_t < x]$ and the Wiener-Hopf factors derived in [13, 12, 11] for Regular Lévy processes of Exponential type (RLPE). Two definitions of RLPE are given: in terms of Lévy densities $f_\pm(x)$ of positive and negative jumps, and in terms of the characteristic exponent $\psi$ of the process. The first definition states that

_S.B._: Department of Economics, The University of Texas at Austin, 2225 Speedway Stop C3100, Austin, TX 78712–0301, sboyarch@utexas.edu
_S.L._: Calico Science Consulting. Austin, TX. Email address: levendorskii@gmail.com.

$f_\pm(x)$ decay exponentially at infinity and obey the asymptotics $f_\pm(x) \sim a_\pm |x|^{-\nu-1}$ as $\pm x \downarrow 0$, where $a_\pm > 0$. According to the second, equivalent, definition, $\psi$ admits analytic continuation to a strip around the real axis, and behaves as $c_\pm |\xi|^\nu$ as $\xi \to \infty$ remaining in the strip and $\pm \operatorname{Re} \xi \to +\infty$ (if $\nu = 0$, as $c_\pm \ln(\pm\xi)$), where $c_\pm > 0$. In the present paper, the asymptotic formulas for the lower tail probability problem are derived under additional conditions on $\psi$ of an RLPE $X$, which are satisfied for wide classes of SINH-regular processes introduced in [14]. The main additional conditions are: $\psi$ admits analytic continuation to a union of a strip and cone around the real axis, and exhibits a regular behavior as $\xi \to \infty$ remaining in the cone. Naturally, the definition of this kind is difficult to reformulate in terms of the Lévy measure. In [18], we introduced the class of Stieltjes-Lévy processes (SL-processes). If $X$ is an SL-process, both the Lévy density and characteristic exponent of an SL process $X$ are defined in terms of certain positive constants (drift can be non-negative) and a pair of Stieltjes-Lévy measures $\mathcal{G}_\pm$ on $\mathbb{R}_+$ (we reproduce the definitions of SINH-regular and SL-processes in Sections 2 and 3). If the measures are signed, the process is called a signed SL-process (sSL process). $\mathcal{G}_\pm$ satisfy certain regularity conditions and $\operatorname{supp} \mathcal{G}_\pm \not\ni 0$ (regular SL-processes), then $X$ is SINH-regular. We say that $X$ is a regular SL process. Stable Lévy processes and mixtures of stable Lévy processes and sSL-processes can be obtained dropping the condition $\operatorname{supp} \mathcal{G}_\pm \not\ni 0$. SL-processes but not sSL processes enjoy an additional property, which simplifies the proofs and the forms of the asymptotic coefficients, namely, for $q > 0$, the equation $q + \psi(\xi) = 0$ has no zeros outside the imaginary axis $i\mathbb{R}$. In [18], we proved that essentially all popular classes of Lévy processes are regular SL-processes, the notable exceptions being stable Lévy processes - a stable Lévy process is an SL process satisfying a weaker form of the regularity condition but a strip of analyticity of $\psi$ is empty, Meixner processes, which are regular sSL-processes but not SL-processes and the Merton model - the process is SINH-regular but not an sSL process. The asymptotic formulas for the survival function are derived for any Lévy process with exponentially decaying tails and refined for SL-processes.

The rest of the paper is organized as follows. In Sections 2 and 3, we recall the definitions of SINH-regular and (s)SL processes, with examples, and several key formulas. In Section 4, we list the extant representations for the Wiener-Hopf factors needed to calculate the asymptotics of $\mathbb{P}[\bar{X}_t < x]$ as $x \to 0$. The asymptotic formulas are derived in Section 5. The asymptotics of the survival probability as $t \to \infty$ is studied in Section 6. In Section 7, we summarize the results and outline possible extensions.

## 2. SINH-REGULAR LÉVY PROCESSES

Let $X$ be a one-dimensional Lévy process on the filtered probability space $(\Omega, \mathcal{F}, \{\mathcal{F}_t\}_{t \geq 0}, \mathbb{P})$ satisfying the usual conditions, and let $\mathbb{E}$ be the expectation operator under $\mathbb{P}$.

2.1. **Definitions.** We use the definition in [10, 12] of the characteristic exponent $\psi(\xi) = \psi^{\mathbb{P}}(\xi)$ of a Lévy process $X$ under $\mathbb{P}$, which is marginally different from the definition in [42]. Namely, $\psi$ is definable from $\mathbb{E}[e^{i\xi X_t}] = e^{-t\psi(\xi)}$. The Lévy-Khintchine formula is

$$(2.1) \qquad \psi(\xi) = \frac{\sigma^2}{2}\xi^2 - i\mu\xi + \int_{\mathbb{R}\setminus 0} (1 - e^{ix\xi} + \mathbf{1}_{(-1,1)}(x) ix\xi) F(dx),$$

where $\sigma^2 \geq 0, \mu \in \mathbb{R}$ and the measure $F(dx)$ satisfies $\int (x^2 \wedge 1) F(dx) < \infty$. $X$ is of finite variation iff $\int (|x| \wedge 1) F(dx) < \infty$, and then the term $\mathbf{1}_{(-1,1)}(x) ix\xi$ can (and will) be omitted.

We represent $\psi$ in the form

$$(2.2) \qquad \psi(\xi) = -i\mu\xi + \psi^0(\xi),$$

where $\mu \in \mathbb{R}$, and impose conditions on $\psi^0$.

We reproduce the definition of SINH-regular processes given in [18, Definition 2.1]; the definition is marginally different from the initial definition in [14]. Also, we impose more stringent conditions on the domain of analyticity assuming that $\mathbb{R}$ is in the interior of the domain of analyticity. Below, $\mu_\pm, \gamma_\pm, \gamma'_\pm, \gamma$ are reals satisfying $\mu_- < 0 < \mu_+$, $-\pi/2 \leq \gamma_- < 0 < \gamma_+ \leq \pi/2$, $\gamma_- < \gamma'_- < \gamma'_+ < \gamma_+$, $\gamma \in (0, \pi]$, $S_{(\mu_-, \mu_+)} = \{\xi \mid \operatorname{Im}\xi \in (\mu_-, \mu_+)\}$ is a strip, and $\mathcal{C}_{\gamma_-, \gamma_+} = \{e^{i\varphi}\rho \mid \rho > 0, \varphi \in (\gamma_-, \gamma_+) \cup (\pi - \gamma_+, \pi - \gamma_-)\}$, $\mathcal{C}_\gamma = \{e^{i\varphi}\rho \mid \rho > 0, \varphi \in (-\gamma, \gamma)\}$ cones.

**Definition 2.1.** *Let $\nu \in (0, 2]$. We say that $X$ is a SINH-regular Lévy process (on $\mathbb{R}$) of type $((\mu_-, \mu_+); \mathcal{C}; \mathcal{C}_+)$ and order $\nu$ iff the following conditions are satisfied:*

*(i) $\mu_- < 0 < \mu_+$;*

*(ii) $\mathcal{C} = \mathcal{C}_{\gamma_-, \gamma_+}, \mathcal{C}_+ = \mathcal{C}_{\gamma'_-, \gamma'_+}$, where $\gamma_- < \gamma'_- < 0 < \gamma'_+ < \gamma_+$;*

*(iii) $\psi^0$ admits analytic continuation to $i(\mu_-, \mu_+) + (\mathcal{C} \cup \{0\})$;*

*(iv) for any $\varphi \in (\gamma_-, \gamma_+)$, there exists $c_\infty(\varphi) \in \mathbb{C} \setminus (-\infty, 0]$ s.t.*

$$(2.3) \qquad \psi^0(\rho e^{i\varphi}) \sim c_\infty(\varphi)\rho^\nu, \quad \rho \to +\infty;$$

*(v) the function $(\gamma_-, \gamma_+) \ni \varphi \mapsto c_\infty(\varphi) \in \mathbb{C}$ is continuous;*

*(vi) for any $\varphi \in (\gamma'_-, \gamma'_+)$, $\operatorname{Re} c_\infty(\varphi) > 0$.*

*We say that $X$ is a SINH-regular Lévy process (on $\mathbb{R}$) of type $((\mu_-, \mu_+); \mathcal{C}; \mathcal{C}_+)$ and order $\nu = 1+$, iff the conditions above bar (2.3) are satisfied, and (2.3) is replaced with*

$$(2.4) \qquad \psi^0(\rho e^{i\varphi}) \sim c_\infty(\varphi)\rho \ln\rho, \quad \rho \to +\infty.$$

*We say that $X$ is a SINH-regular Lévy process (on $\mathbb{R}$) of type $((\mu_-, \mu_+); \mathcal{C}; \mathcal{C}_+)$ and order $\nu = 0+$, iff the conditions (i)-(iii) are satisfied, and there exists $c > 0$ such that*

$$(2.5) \qquad \psi^0(\xi) \sim c \ln|\xi|, \quad (i(\mu_-, \mu_+) + \mathcal{C} \ni)\xi \to \infty.$$

**Remark 2.1.** Conditions for $\varphi \in (\pi - \gamma_-, \pi - \gamma_+)$ and $\varphi \in (\pi - \gamma'_-, \pi - \gamma'_+)$ follow from the conditions for $\varphi \in (\gamma_-, \gamma_+)$ and $\varphi \in (\gamma'_-, \gamma'_+)$ because $\psi(-\bar{\xi}) = \overline{\psi(\xi)}$. In Definition 2.1, conditions are imposed on $\psi^0$ whereas in [14], the same conditions are imposed on $\psi$. If either $\nu \in (1, 2]$ or $\mu = 0$, then the conditions on $\psi^0$ and $\psi$ are equivalent, and the process is an elliptic[1] SINH-process of order $\nu$ in the terminology of [14]. If $\nu < 1$ and $\mu \neq 1$, then the process $X$ with the characteristic exponent $-i\mu\xi + \psi^0(\xi)$ is elliptic of order 1, in the terminology of [14]. Furthermore, according to the definition in [14], $\mathcal{C}_+$ is determined by the drift term if $\mu \neq 0$: $\mathcal{C}_+$ is the intersection of $\mathcal{C}$ with the upper (resp. lower) half-plane if $\mu > 0$ (resp. $\mu < 0$).

Finally, note that in [18], the reader can find the definition of a more general class of SINH-regular processes, with examples.

2.2. **Examples of SINH-regular processes.** Essentially all popular Lévy processes are SINH-regular. We reproduce the list from [18, Sect. 3.3].

---

[1]The name elliptic is natural from the point of view of the theory of PDO: if (2.3) holds, then the infinitesimal generator $L^0 = -\psi^0(D)$ is an elliptic PDO.

2.2.1. *The Brownian motion (BM)*. BM is of order 2; since $\psi^0(\xi) = \frac{\sigma^2}{2}\xi^2$ is an entire function, $\mathcal{C} = \mathbb{C}$, $\mu_- = -\infty$, $\mu_+ = +\infty$. For any $\varphi \in [0, 2\pi)$,

$$(2.6) \qquad \psi^0(\rho e^{i\varphi}) \sim \frac{\sigma^2}{2}\rho^2 e^{2i\varphi}, \ \rho \to +\infty,$$

hence, $c_\infty(\varphi) = \frac{\sigma^2}{2}e^{2i\varphi}$, and $\operatorname{Re} c_\infty(\varphi) > 0 \Leftrightarrow \cos(2\varphi) > 0$. It follows that $\mathcal{C}_+ = \mathcal{C}_{-\pi/4, \pi/4}$.

2.2.2. *Merton model* [40]. The characteristic exponent is given by

$$(2.7) \qquad \psi^0(\xi) = \frac{\sigma^2}{2}\xi^2 + \lambda \cdot \left(1 - e^{im\xi - \frac{s^2}{2}\xi^2}\right),$$

where $\sigma, s, \lambda > 0$ and $\mu, m \in \mathbb{R}$. As far as the analytical properties formulated in the definition of SINH-processes are concerned, the difference with BM is that $\mathcal{C} = \mathcal{C}_+ = \mathcal{C}_{-\pi/4, \pi/4}$, and $\mathcal{C} = \mathcal{C}_{\gamma_-, \gamma_+}$ with either $\gamma_- < -\pi/4$ or $\gamma_+ > \pi/4$ cannot be used.

2.2.3. *Lévy processes with rational characteristic exponents and non-trivial BM component.* The order is 2, and an admissible strip of analyticity $S_{(\mu_-, \mu_+)}$ around the real axis may not contain poles of $\psi^0$. Explicit formulas for the Wiener-Hopf factors are easy to derive (see, e.g., [12]) in terms of the poles of $\psi^0$ and zeros of the function $q + \psi(\xi)$, the multiplicities of zeros and poles being taken into account. After $\mu_-, \mu_+$ are chosen, $\mathcal{C}$ is the maximal cone around the real axis such that the region $i(\mu_-, \mu_+) + (\mathcal{C} \cup \{0\})$ contains no poles, and $\mathcal{C}_+ = \mathcal{C} \cap \mathcal{C}_{-\pi/4, \pi/4}$. Since the efficiency of SINH-acceleration depends, mostly, on the opening angle of $\mathcal{C}_+$, it is advisable to choose small (in absolute value) $\mu_-$ and $\mu_+$ so that $\mathcal{C}_+$ can be chosen "wider". Lévy processes of the phase type [1, 2] have rational characteristic exponents, hence, the recommendations above are applicable.

Calculation of the rational characteristic exponent is straightforward if the Lévy densities of positive and negative jumps are mixtures of exponential polynomials. Furthermore, all poles are on $i\mathbb{R}$, hence, $\mathcal{C} = \mathbb{C} \setminus i\mathbb{R}$ [32, 33]. The factorization of $q + \psi(\xi)$ (calculation of the Wiener-Hopf factors) simplifies if all the roots of the characteristic equation $\psi(\xi) + q = 0$ are on the imaginary axis. Then the roots can be easily calculated, and explicit formulas for the Wiener-Hopf factors as sums or products derived. See, e.g., [32, 33]. A popular special case is the hyper-exponential jump-diffusion model (HEJD model) introduced in [32] without a special name assigned and in [37], and studied in detail in [32, 33]). The Lévy measure is of the form

$$(2.8) \qquad F(dx) = \mathbf{1}_{(-\infty, 0)}(x)\sum_{k=1}^{n^-} p_j^- \alpha_j^- e^{\alpha_j^- x} + \mathbf{1}_{(0, +\infty)}(x)\sum_{j=1}^{n^+} p_j^+ \alpha_j^+ e^{-\alpha_j^+ x},$$

where $n^\pm$ are positive integers, and $\alpha_j^\pm, p_j^\pm > 0$ are reals. The characteristic exponent is

$$(2.9) \qquad \psi^0(\xi) = \frac{\sigma^2}{2}\xi^2 - i\mu\xi + \sum_{j=1}^{n^+} p_j^+ \frac{-i\xi}{\alpha_j^+ - i\xi} + \sum_{k=1}^{n^-} p_k^- \frac{i\xi}{\alpha_k^- + i\xi}.$$

Double-exponential jump diffusion model introduced to finance in [28] (and well-known for decades) is a special case of HEJD models with $n^+ = n^- = 1$. The order is 2, $\mu_+ = \min \alpha_k^-$, $\mu_- = -\min \alpha_k^+$, and $\mathcal{C}, \mathcal{C}_+$ are as in the BM model.

In [24], a class of processes with the Lévy measure of the form (2.8) with some of $p_j^\pm$ being negative is introduced, and the name mixed exponential jump diffusion model (MEJD) is suggested. Sufficient conditions for $p_j^{\pm}$ and $\alpha_j^{\pm}$ to define the non-negativity of the densities are $p_1^{\pm} > 0$ and $\sum_{j=1}^{k} p_j^{\pm} \alpha_j^{\pm} \geq 0$, $k = 1, 2, \ldots, n^{\pm}$. An important qualitative difference between HEJD and MEJD is that in HEJD models, the Lévy densities of positive and negative jumps are monotone (in fact, completely monotone), whereas in MEJD, the densities may be non-monotone. Note that the Lévy densities given by mixtures of exponential polynomials [32, 33] are typically non-monotone, and qualitative properties of MEJD densities can be easily reproduced by exponential polynomials. As the simplest example, the reader can compare the following two functions on $\mathbb{R}_+$: $f_1(x) = e^{-\lambda_1 x} - e^{-\lambda_2 x}$, where $0 < \lambda_1 < \lambda_2$, and $f_2(x) = xe^{-\lambda_1 x}$.

2.2.4. *Variance Gamma processes (VGP)*. VG model was introduced to Finance in [38]. The characteristic exponent can be written in the form

$$(2.10) \qquad \psi^0(\xi) = c[\ln(\alpha^2 - (\beta + i\xi)^2) - \ln(\alpha^2 - \beta^2)],$$

where $\alpha > |\beta| \geq 0$, $c > 0$. VGP is SINH-regular of type $((-\alpha + \beta, \alpha + \beta); \mathbb{C} \setminus i\mathbb{R}, \mathbb{C} \setminus i\mathbb{R})$ and order 0+ because $\forall \; \varphi \in (-\pi/2, \pi/2)$,

$$(2.11) \qquad \psi^0(\rho e^{i\varphi}) = c(\ln \rho + i\varphi) + O(1), \;\; \rho \to +\infty.$$

2.2.5. *NIG and NTS*. Normal inverse Gaussian (NIG) processes, and the generalization: Normal Tempered Stable (NTS) processes are constructed in [4, 5], respectively. The characteristic exponent is given by

$$(2.12) \qquad \psi^0(\xi) = \delta[(\alpha^2 - (\beta + i\xi)^2)^{\nu/2} - (\alpha^2 - \beta^2)^{\nu/2}],$$

where $\nu \in (0, 2)$, $\delta > 0$, $|\beta| < \alpha$; NIG obtains with $\nu = 1$. This is a SINH-regular process of order $\nu$ and type $((-\alpha + \beta, \alpha + \beta); \mathbb{C} \setminus i\mathbb{R}, \mathcal{C}_{-\gamma_\nu, \gamma_\nu})$, where $\gamma_\nu = \min\{1, 1/\nu\}\pi/2$. Indeed, for $\varphi \in (-\pi/2, \pi/2)$,

$$(2.13) \qquad \psi^0(\rho e^{i\varphi}) = \delta e^{i\varphi\nu}\rho^\nu + O(\rho^{\nu-1}) + O(1), \;\; \rho \to +\infty.$$

2.2.6. *The Meixner process*. For the background, see, e.g., [44, 41, 39]. The Lévy density of the Meixner process $X$ is

$$(2.14) \qquad f(x) = \delta \frac{\exp(bx/a)}{x\sinh(\pi x/a)},$$

where $\delta, a > 0$ and the asymmetry parameter $b \in (-\pi, \pi)$. The characteristic exponent is

$$(2.15) \qquad \psi^0(\xi) = 2\delta[\ln[\cosh((a\xi - ib)/2))] - \ln\cos(b/2)].$$

The formula $\ln \cosh(z) = z + \ln(1 + e^{-2z}) - \ln 2$ defines a function analytic in the right half-plane, hence, $\psi^0(\xi)$ admits analytic continuation to $\mathbb{C} \setminus i\mathbb{R}$. Set $\mu_- = (-\pi + b)/a$, $\mu_+ = (\pi + b)/a$. Since $\cosh((a\xi - ib)/2)) > 0$ for $\xi \in i(\mu_-, \mu_+)$, $\psi^0$ is analytic in $\mathbb{C} \setminus i((-\infty, \mu_-] \cup [\mu_+, +\infty))$. Let $\gamma \in (0, \pi/2)$. As $\xi \to \infty$ in $\mathcal{C}_\gamma$, $\psi(\xi) \sim a\delta\xi$, therefore, $X$ is SINH-regular of order 1 and type $(((-\pi + b)/a, (\pi + b)/a), \mathbb{C} \setminus i\mathbb{R}, \mathbb{C} \setminus i\mathbb{R})$.

2.2.7. *KoBoL processes.* A generic process of Koponen's family [9, 10] was constructed as a mixture of spectrally negative and positive pure jump processes, with the Lévy measure

$$(2.16) \qquad F(dx) = c_+ e^{\lambda_- x} x^{-\nu_+ -1} \mathbf{1}_{(0,+\infty)}(x) dx + c_- e^{\lambda_+ x} |x|^{-\nu_- -1} \mathbf{1}_{(-\infty,0)}(x) dx,$$

where $c_\pm > 0, \nu_\pm \in [0,2), \lambda_- < 0 < \lambda_+$. In this paper, we allow for $c_+ = 0$ or $c_- = 0$, $\lambda_- = 0 < \lambda_+$ and $\lambda_- < 0 \le \lambda_+$. This generalization is almost immaterial for evaluation of probability distributions and expectations because for efficient calculations, the first crucial property, namely, the existence of a strip of analyticity of the characteristic exponent, around or adjacent to the real line, holds if $\lambda_- < \lambda_+$ and $\lambda_- \le 0 \le \lambda_+$. [2] Furthermore, the Esscher transform allows one to reduce both cases $\lambda_- = 0 < \lambda_+$ and $\lambda_- < 0 \le \lambda_+$ to the case $\lambda_- < 0 < \lambda_+$. Using the Lévy-Khintchine formula, it is straightforward to derive from (2.16) explicit formulas for $\psi^0$. See [9, 10, 12]. If $\nu_\pm \in (0,2), \nu_\pm \ne 1$,

$$(2.17) \qquad \psi^0(\xi) = c_+ \Gamma(-\nu_+)((-\lambda_-)^{\nu_+} - (-\lambda_- - i\xi)^{\nu_+}) + c_- \Gamma(-\nu_-)(\lambda_+^{\nu_-} - (\lambda_+ + i\xi)^{\nu_-}).$$

(Formulas in the case $\nu_\pm = 0, 1$ are in [18].) A subclass with $\nu_+ = \nu_- = \nu \in (0,2)$ and $c_+ = c_-$ was labelled KoBoL in [12] and called a process of order $\nu$. To simplify the name, we will call a pure jump process with the Lévy measure (2.16) a KoBoL process as well. As in [10], we allow for $\nu_- \ne \nu_+$ and $c_+ \ne c_-$. If either $c_- = 0$ or $c_+ = 0$, we say that the process is a one-sided KoBoL. The formula for one-sided KoBoL of order 1 is different. See [18]. One-sided KoBoL processes were used in [8] to price CDSs. Mixing one-sided processes of order $\nu \in (0,2), \nu \ne 1$, and order 1, one can obtain more exotic characteristic exponents.

Note that a specialization $\nu \ne 1, c = c_\pm > 0$, of KoBoL used in a series of numerical examples in [9, 12] was named CGMY model in [25] (and the labels were changed: letters $C, G, M, Y$ replace the parameters $c, \nu, \lambda_-, \lambda_+$ of KoBoL):

$$(2.18) \qquad \psi^0(\xi) = c\Gamma(-\nu)[(-\lambda_-)^\nu - (-\lambda_- - i\xi)^\nu + \lambda_+^\nu - (\lambda_+ + i\xi)^\nu].$$

Evidently, $\psi^0$ given by (2.18) is analytic in $\mathbb{C} \setminus i\mathbb{R}$, and $\forall \varphi \in (-\pi/2, \pi/2)$, (2.3) holds with

$$(2.19) \qquad c_\infty(\varphi) = -2c\Gamma(-\nu)\cos(\nu\pi/2)e^{i\nu\varphi}.$$

Hence, $X$ is SINH-regular of type $((\lambda_-, \lambda_+), \mathbb{C} \setminus \{0\}, \mathcal{C}_{-\gamma_\nu, \gamma_\nu})$, where $\gamma_\nu = \min\{1, 1/\nu\}\pi/2$, and order $\nu$. For the calculation of order and type for a generic KoBoL, see [18].

2.2.8. *The $\beta$-class* [29]. The characteristic exponent is of the form

$$(2.20) \qquad \begin{aligned} \psi^0(\xi) & = \frac{\sigma^2}{2}\xi^2 + \frac{c_1}{\beta_1}\left\{B(\alpha_1, 1-\gamma_1) - B\left(\alpha_1 - \frac{i\xi}{\beta_1}, 1-\gamma_1\right)\right\} \\ & \quad + \frac{c_2}{\beta_2}\left\{B(\alpha_2, 1-\gamma_2) - B\left(\alpha_2 + \frac{i\xi}{\beta_2}, 1-\gamma_2\right)\right\}, \end{aligned}$$

where $c_j \ge 0, \alpha_j, \beta_j > 0$ and $\gamma_j \in (0,3) \setminus \{1,2\}$, and $B(x,y)$ is the Beta-function. It was shown in [29], that all poles of $\psi^0$ are on $i\mathbb{R}_+ \setminus 0$. Hence, $(\mu_-, \mu_+)$ is the maximal interval containing 0 and no poles of $\psi^0$, and $\mathcal{C} = \mathbb{C} \setminus i\mathbb{R}$. For calculation of the order of the process and $\mathcal{C}_+$ as functions of the parameters in (2.20). The number of parameters is larger than in the case of KoBoL but the variety of different cases (order and type) is essentially the same as in the case of KoBoL.

---

[2] The property does not hold if there is no such a strip (formally, $\lambda_- = 0 = \lambda_+$). The classical example are stable Lévy processes. The conformal deformation technique can be modified for this case [15, 17].

2.2.9. *Meromorphic Lévy processes* [30]. The Lévy measures and characteristic exponents of the meromorphic Lévy processes introduced in [30] are defined by almost the same formulas as in HEJD model. The difference is that the sums in (2.8) and (2.9) are infinite. A natural condition $\alpha_j^\pm \to +\infty$ as $j \to +\infty$ and the requirement that the infinite sum defines a Lévy measure $\sum_{j\geq 0} p_j^\pm (\alpha_j^\pm)^{-2} < +\infty$ are imposed. The poles of $\psi(\xi)$ are on $\mathbb{R} \setminus 0$. If $\sigma^2 > 0$, meromorphic processes are SINH-regular of order 2 and type $((\mu_-\mu_+), \mathbb{C} \setminus i\mathbb{R}, \mathcal{C}_{-\pi/4,\pi/4})$, where $(\mu_-,\mu_+)$ is the maximal interval containing 0 and no poles of $\psi^0$. If $\sigma^2 = 0$, additional conditions on the asymptotics of $p_j^\pm$ and $\alpha_j^\pm$ as $j \to +\infty$ need to be imposed to obtain a SINH-regular process.

## 3. Stieltjes-Lévy processes

We reproduce several definitions, examples and results from [18].

3.1. **Preliminaries.** Let the Lévy density be absolutely continuous: $F(dx) = f(x)dx$, and let $\psi$ admit analytic continuation to $\mathbb{C} \setminus i\mathbb{R}$. Then, under additional conditions, we can express $f_+ = \mathbf{1}_{(0,+\infty)}f$ and $f_- = \mathbf{1}_{(-\infty,0)}f$ in terms of integrals over the cuts $i[(-\infty,\mu_-]$ and $i[\mu_+,+\infty)$, respectively, w.r.t. to certain measures $\mathcal{G}_\pm(dt) = \mathcal{G}_\pm(\psi; dt)$ (possibly, *signed*):

$$(3.1) \qquad f_+(x) \;=\; \int_{(0,+\infty)} e^{-tx} \mathcal{G}_+(dt), \; x > 0,$$

$$(3.2) \qquad f_-(x) \;=\; \int_{(0,+\infty)} e^{tx} \mathcal{G}_-(dt), \; x < 0,$$

where $\text{supp}\,\mathcal{G}_+ \subset [-\mu_-, +\infty)$ (if $\mu_- = 0$, $\mathcal{G}_+$ has no atom at 0), and $\text{supp}\,\mathcal{G}_- \subset [\mu_+, +\infty)$ (if $\mu_+ = 0$, $\mathcal{G}_-$ has no atom at 0). In terms of the Laplace transforms $\widetilde{\mathcal{G}_\pm(dt)}$ of $\mathcal{G}_\pm(dt)$,

$$(3.3) \qquad f_+(x) = \widetilde{\mathcal{G}_+(dt)}(x), \quad f_-(x) = \widetilde{\mathcal{G}_-(dt)}(-x).$$

**Example 3.1.** In the case of HEJD, the $\beta$-model and meromorphic processes in general, the measures $\mathcal{G}_\pm(dt) = \mathcal{G}_\pm(\psi; dt)$ are discrete:

$$(3.4) \qquad \mathcal{G}_\pm(dt) = \sum_{\alpha \in \mathcal{A}_\pm} g_\alpha^\pm \delta_\alpha,$$

where $\mathcal{A}_\pm$ are finite or discrete sets with the only accumulation point at $+\infty$. The set $-\mathcal{A}_+$ (resp., $\mathcal{A}_-$) is the set of (simple) poles of $\psi$ on $(-\infty, 0)$ (resp., $(0, +\infty)$), and $g_\alpha^+ = \text{Res}(i\psi, -i\alpha)$, $g_\alpha^- = -\text{Res}(i\psi, i\alpha)$ are positive. If we relax the conditions on the parameters of HEJD and meromorphic model so that some of $g_\alpha^\pm$ are negative (but the process is a Lévy process), then we obtain the representations (3.1)-(3.2) with signed measures.

**Example 3.2.** Let $X$ be the one-sided stable Lévy process of index $\alpha \in (0,2), \alpha \neq 1$, with the Lévy density $f_+(x) = x^{-\alpha-1}\mathbf{1}_{x>0}$ and the characteristic exponent $\psi_+(\xi) = -\Gamma(-\alpha)(0 - i\xi)^\alpha$. Then (3.1) holds with $\mathcal{G}_+(dt) = \Gamma(-\alpha)\sin(-\pi\alpha)\pi^{-1}t^\alpha dt$.

**Definition 3.3.** $X^\pm = X_{\mathcal{G}_\pm}^\pm$ *denote the one-sided Lévy processes given by the generating triplets* $(0, 0, f_\pm(x)dx)$, *where* $f_\pm(x) := f_\pm(\mathcal{G}_\pm; x)$ *are defined by (3.1)-(3.2). The characteristic exponents of* $X_{\mathcal{G}_\pm}^\pm$ *are denoted* $\psi_\pm(\xi) := \psi_\pm(\mathcal{G}_\pm, \xi)$.

Evidently, $\psi_-(\mathcal{G}; \xi) = \psi_+(\mathcal{G}, -\xi)$, therefore, it suffices to derive the condition for $\mathcal{G}$ to define the Lévy densities for $X_{\mathcal{G}}^+$ and formula for $\psi_+(\mathcal{G}; \xi)$; the condition for $X_{\mathcal{G}}^-$ is the same, and the formula for $\psi_-(\mathcal{G}; \xi)$ obtains by symmetry.

**Lemma 3.4.** *(a) Let $\mathcal{G}_+ \geq 0$. Then (3.1) defines a Lévy density if and only if*

$$(3.5) \qquad \int_{(0,+\infty)} \frac{\mathcal{G}_+(dt)}{t + t^m} < \infty,$$

*where $m = 3$; the density is completely monotone.*
*(b) If*

$$(3.6) \qquad \int_{(0,+\infty)} \frac{|\mathcal{G}(dt)|}{t + t^m} < \infty,$$

*where $m = 3$, and $f_+$ given by (3.1) is non-negative, then $f_+$ is a Lévy density.*
*(c) The pure jump process $X^+$ with the Lévy density $f_+$ is of finite variation iff (3.6) holds with $m = 2$.*

### 3.2. Stieltjes-Lévy measures and functions. Stieltjes transform.

**Definition 3.5.** *A non-negative measure $\mathcal{G}$ on $(0, +\infty)$ is a Stieltjes measure iff*

$$(3.7) \qquad \int_{(0,+\infty)} (1+t)^{-1} \mathcal{G}(dt) < \infty.$$

*We write $\mathcal{G} \in \mathrm{SM}_0$. The Stieltjes transform $ST(\mathcal{G})$ of $\mathcal{G}$ is given by*

$$(3.8) \qquad ST(\mathcal{G})(z) = \int_{(0,+\infty)} (z+t)^{-1} \mathcal{G}(dt).$$

The definitions of the Stieltjes measure and transform are standard - see, e.g., ([43, Defin. 2.1]). We introduce the notation $ST(\mathcal{G})$ to shorten the definitions, statements and proofs.

**Definition 3.6.** *If $\mathcal{G} \in \mathrm{SM}_0$ and $\mathrm{supp}\,\mathcal{G} \subset [\mu, +\infty)$, where $\mu > 0$, we write $\mathcal{G} \in \mathrm{SM}_\mu$.*
*If $\mathcal{G} \in \mathrm{SM}_\mu$, $\mu \geq 0$, and $f = ST(\mathcal{G}_\mu)$, we write $f \in \mathcal{S}_\mu$.*

Evidently, for any $\mu' \in [0, \mu]$, $\mathrm{SM}_\mu \subset \mathrm{SM}_{\mu'}$, $\mathcal{S}_{\mu'} \subset \mathcal{S}_\mu$.

**Proposition 3.7.** *Let $\mathcal{G} \in \mathrm{SM}_\mu$. Then $ST(\mathcal{G})$ is analytic in $\mathbb{C} \setminus (-\infty, -\mu]$.*

**Proposition 3.8.** *Measure $\mathcal{G}$ satisfies (3.5) with $m = 3$ (resp., with $m = 2$) iff there exist $a_2, a_1 \geq 0$ and $\mathcal{G}^0 \in \mathrm{SM}_0$ such that $\mathcal{G}(dt) = (a_2 t^2 + a_1 t)\mathcal{G}^0(dt)$ (resp., $\mathcal{G}(dt) = t\mathcal{G}^0(dt)$).*

### 3.3. Definition of SL and sSL processes.
The class $\mathcal{S}$ of Stieltjes functions (see, e.g., [43, Defin. 2.1]) is wider than $\mathcal{S}_0$: $f \in \mathcal{S}$ if there exist $a_0, a_1 \geq 0$ and $\mathcal{G} \in \mathrm{SM}_0$ such that $f(z) = a_0/z + a_1 + ST(\mathcal{G})(z)$. For construction of spectrally one-sided Lévy processes, $\mathcal{S}$ is appropriate. Indeed, any $\mathcal{G} \in \mathrm{SM}_0$ is the Stieltjes measure of a complete Bernstein function $g$, which appears in the Stieltjes representation of $g$: $g(z) = a_0 + a_1 z + z ST(\mathcal{G})(z)$ (see, e.g., [43, Thm 6.2, Corr. 6.3 and Remark 6.4]). Evidently, if $a_0 = 0$, $\psi(\xi) = g(-i\xi)$ is the characteristic exponent of a subordinator. For construction of more general Lévy processes, class $\mathcal{S}_0$ is more convenient.

**Definition 3.9.** *Let $\mu \geq 0$. We say that the measure $\mathcal{G}$ on $[\mu, +\infty)$ is a Stieltjes-Lévy measure (SL measure) of class $\mathrm{SLM}_\mu$ if there exists $a_2, a_1 \geq 0$, $a_2 + a_1 > 0$, such that $(a_2 t^2 + a_1 t)^{-1} \mathcal{G}(dt) \in \mathrm{SM}_\mu$.*

*We say that $\mathcal{G}$ is a signed Stieltjes-Lévy measure (sSL measure) of class $\mathrm{sSLM}_\mu$ if $\mathcal{G}$ admits the representation $\mathcal{G} = \mathcal{G}_1 - \mathcal{G}_2$, where $\mathcal{G}_j \in \mathrm{SLM}_\mu$, and the Laplace transform $\tilde{\mathcal{G}}$ of $\mathcal{G}$ is non-negative on $(0, +\infty)$.*

**Definition 3.10.** *Let $\mathcal{G}^0 \in \mathrm{SM}_\mu$, $\mathcal{G}(dt) = (a_2 t^2 + a_1 t)\mathcal{G}^0(dt)$, and $X^\pm = X^\pm_{\mathcal{G}}$. Then we write $X^\pm \in SL^\pm_\mu$. If $a_1 = 0$ (resp., $a_2 = 0$), we write $X^\pm \in SL^{2;\pm}_\mu$ (resp., $X^\pm \in SL^{1;\pm}_\mu$).*

Evidently, if $\mu > 0$, one can use a simpler definition of classes $SL^{2;\pm}_\mu$ instead of the general definition of classes $SL^\pm_\mu$; the statement $X \in SL^{1;\pm}_\mu$ is useful if we wish to indicate that the jump component of $X$ is of finite variation. The following proposition demonstrates that $SL^\pm_0 \neq SL^{1;\pm}_0 \cup SL^{2;\pm}_0$.

**Proposition 3.11.** *Let $X^\pm$ be the one-sided stable Lévy process of index $\alpha \in (0, 2)$, with the Lévy density $f_\pm(x) = \Gamma(\alpha + 1)|x|^{-\alpha-1}, \pm x > 0$. Then*
*(a) if $\alpha \in (0, 1)$, then $X^\pm$ is in $SL^{1;\pm}_0$ but not in $SL^{2;\pm}_0$;*
*(b) if $\alpha \in (1, 2)$, then $X^\pm$ is in $SL^{2;\pm}_0$ but not in $SL^{1;\pm}_0$;*
*(c) if $\alpha = 1$, then $X^\pm$ is in $SL^\pm_0$ but not in $SL^{1;\pm}_0 \cup SL^{2;\pm}_0$.*

In [18], we proved that stable Lévy processes, KoBoL, VGP, NIG, and NTS are SL-processes, whereas Meixner processes are sSL processes, and derive explicit formulas for the corresponding SL- and sSL-measures.

**Theorem 3.12.** *Let $\mathcal{G}(dt) = (a_2 t^2 + a_1 t)\mathcal{G}^0(dt) \in SLM_\mu$ and $X^\pm = X^\pm_{\mathcal{G}}$. Then*
*(a) the characteristic exponent of $X^\pm$ is of the form*

$$(3.9) \qquad \psi_\pm(\xi) = (a_2 \xi^2 \mp i a_1 \xi) ST(\mathcal{G}^0)(\mp i\xi) \pm ic\xi,$$

*where $c = c(a_2, a_1, \mathcal{G}^0) \in \mathbb{R}$. If $a_2 = 0$, and the Lévy-Khintchine formula for processes of finite variation is used, then $c = 0$;*
*(b) $X^\pm \in SL^{1;\pm}_\mu$ are finite variation processes;*
*(c) if $a_2 = 0$ and $t\mathcal{G}^0(dt) \in L_1$, then $X^\pm \in SL^{1;\pm}_1$ are of finite activity;*
*(d) if $\mu > 0$, then $SL^\pm_\mu = SL^{2;\pm}_\mu$;*
*(e) if $X^\pm \in SL^\pm_\mu$ and $t\mathcal{G}^0(dt) \in \mathrm{SM}_\mu$, then $X^\pm \in SL^{1;\pm}_\mu$.*

In the following Lemma, by a slight abuse of notation, we denote by $\psi_\pm$ the characteristic exponent defined by the generating triplet $(\sigma^2, b, f_\pm(x)dx)$, where $f_\pm(x) = f_\pm(\mathcal{G}; x)$ are defined by (3.1)-(3.2) with $\mathcal{G}_\pm = \mathcal{G}(a_2, a_1, \mathcal{G}^0)$. In the case $a_2 = 0$, we use the Lévy-Khintchine formula for jump component of finite variation.

**Lemma 3.13.** *(a) $\psi_+$ is analytic in $\mathbb{C} \setminus i(-\infty, -\mu]$, $\psi_-$ is analytic in $\mathbb{C} \setminus i[\mu, +\infty)$;*
*(b) $\forall \gamma \in (0, \pi)$, $\psi_\pm(\xi) \sim \frac{\sigma^2}{2} \xi^2$ as $(\pm i\mathcal{C}_\gamma \ni)\xi \to \infty$, uniformly in $\arg\xi \in [-\gamma, \gamma]$.*
*(c) If $\sigma^2 = a_2 = 0$, then, $\forall \gamma \in (0, \pi)$, $\psi_\pm(\xi) \sim -ib\xi$ as $(\pm i\mathcal{C}_\gamma \ni)\xi \to \infty$, uniformly in $\arg\xi \in [-\gamma, \gamma]$.*

**Definition 3.14.** *Let $\mu_- \leq 0 \leq \mu_+$. A Lévy process $X$ is called a signed Stieltjes-Lévy process (sSL-process) of class $sSL_{\mu_-,\mu_+}$ if the Lévy density of $X$ is of the form*

$$(3.10) \qquad f(x) = \mathbf{1}_{(-\infty,0)}(x) \int_{(0,+\infty)} e^{tx} \mathcal{G}_-(dt) + \mathbf{1}_{(0,+\infty)}(x) \int_{(0,+\infty)} e^{-tx} \mathcal{G}_+(dt),$$

*where $\mathcal{G}_- \in \mathrm{sSLM}_{\mu_+}$ and $\mathcal{G}_+ \in \mathrm{sSLM}_{-\mu_-}$. If $\mathcal{G}_- \in \mathrm{SLM}_{\mu_+}$ and $\mathcal{G}_+ \in \mathrm{SLM}_{-\mu_-}$, $X$ is called a Stieltjes-Lévy process (SL-process) of class $SL_{\mu_-,\mu_+}$.*

The following theorem is immediate from Definition 3.14 and Lemma 3.13.

**Theorem 3.15.** *Let $\psi$ be the characteristic exponent of $X \in sSL_{\mu_-,\mu_+}$. Then*
*1. $\psi$ is of the form*

$$(3.11) \qquad \psi(\xi) = (a_2^+ \xi^2 - ia_1^+ \xi) ST(\mathcal{G}_+^0)(-i\xi) + (a_2^- \xi^2 + ia_1^- \xi) ST(\mathcal{G}_-^0)(i\xi) + (\sigma^2/2)\xi^2 - i\mu\xi,$$

*where $\mathcal{G}_+^0 \in SM_{-\mu_-}, \mathcal{G}_+^0 \in SM_{\mu_+}$;*
*2. $\psi$ admits analytic continuation to $\mathbb{C} \setminus i((-\infty, \mu_-] \cup [\mu_+, +\infty))$;*
*3. $\forall\ \gamma \in (0, \pi/2),\ \psi(\xi) \sim \frac{\sigma^2}{2}\xi^2$ as $\xi \to \infty$, uniformly in $\arg \xi \in [-\gamma, \gamma] \cup [\pi - \gamma, \pi + \gamma]$;*
*4. if $\sigma^2 = a_2^+ = a_2^- = 0$, then, for any $\gamma \in (0, \pi/2)$, $\psi(\xi) \sim -ib\xi$ as $\xi \to \infty$, uniformly in $\arg \xi \in [-\gamma, \gamma] \cup [\pi - \gamma, \pi + \gamma]$.*

### 3.4. Representations of sSL-measures in terms of the characteristic exponent.
Let $\mathcal{G}_+ = \mathcal{G}_{+;+} - \mathcal{G}_{+;-}$, $\mathcal{G}_{+;\pm} \in \mathrm{SLM}_{-\mu_-}$ and $\mathcal{G}_- = \mathcal{G}_{-;+} - \mathcal{G}_{-;-}$, $\mathcal{G}_{-;\pm} \in \mathrm{SLM}_{\mu_+}$ be the Jordan decompositions of measures $\mathcal{G}_\pm$ in (3.10) Denote by $U_{+;\pm}$ the set of points of continuity of $t \mapsto \mathcal{G}_{+;\pm}(-\infty, t)$ and set $U_+ = U_{+;+} \cap U_{+,-}$. Similarly, define $U_-$.

**Theorem 3.16.** *Let $X$ be of class $sSL_{\mu_-,\mu_+}$, with the Lévy density (3.10) and characteristic exponent $\psi$. Then*
*(a) for any $\theta \in (\mu_-, \mu_+)$ and $x > 0$,*

$$(3.12) \qquad f_\pm(\pm x) \;=\; \frac{1}{\pi} \lim_{\epsilon \to 0+} \int_{-\theta}^{+\infty} e^{-tx} \operatorname{Im} \psi(\mp(it + \epsilon)) dt;$$

*(b) for any $u, v \in U_\pm$,*

$$(3.13) \qquad \mathcal{G}_\pm((u,v]) = \lim_{\epsilon \to 0+} \frac{1}{\pi} \int_u^v \operatorname{Im} \psi(\mp(it + \epsilon)) dt;$$

*(c) if $\mathcal{G}_\pm(\{\mp\mu_\mp\}) = 0$, then (3.12) holds with $\mp\mu_\mp$ in place of $-\theta$.*

### 3.5. Examples of sSL- and SL-processes.

**Example 3.17.** *Let $\psi = \psi^0$ be given by (2.12) (the case of NIG and NTS processes). Then*

$$(3.14) \qquad \operatorname{Im} \psi(it + 0) \;=\; \delta \sin(\pi\nu/2)((t - \beta)^2 - \alpha^2)^{\nu/2},\ t > \mu_+ := \alpha + \beta,$$

$$(3.15) \qquad \operatorname{Im} \psi(it - 0) \;=\; \delta \sin(\pi\nu/2)((t + \beta)^2 - \alpha^2)^{\nu/2},\ t < \mu_- := -\alpha + \beta.$$

Hence, $X$ is an SL process, which is a mixture of two independent one-sided SL-processes.

**Example 3.18.** Let $\psi = \psi^0$ be given by (2.17) (the case of KoBoL with the positive and negative densities of order $\nu_{\pm} \in (0, 2) \setminus \{1\}$). Then $\psi$ is the characteristic exponent of an SL process, and

$$\tag{3.16} \operatorname{Im} \psi(it + 0) = -c_- \Gamma(-\nu_-) \sin(\pi \nu_-)(t - \lambda_+)^{\nu_-}, \ t > \lambda_+,$$

$$\tag{3.17} \operatorname{Im} \psi(-it - 0) = -c_+ \Gamma(-\nu_+) \sin(\pi \nu_+)(t + \lambda_-)^{\nu_+}, \ t > -\lambda_-.$$

If $\nu_{\pm} = 1$, then

$$\tag{3.18} \operatorname{Im} \psi(it + 0) = \pi c_-(t - \lambda_+), \ t > \lambda_+,$$

$$\tag{3.19} \operatorname{Im} \psi(-it - 0) = \pi c_+(t + \lambda_-), \ t > -\lambda_-.$$

**Example 3.19.** Let $\psi = \psi^0$ be the characteristic exponent of a KoBoL process of order 0+ (the asymmetric version of VGP). Then

$$\tag{3.20} \operatorname{Im} \psi(it \pm 0) = \pi c_{\mp}, \ \pm t > \pm \lambda_{\pm},$$

and $\psi$ is the characteristic exponent of a regular SL-process.

**Example 3.20.** Let $\psi^0$ be given by (2.15). For $\epsilon \neq 0$, and $t \in \mathbb{R}$, we have

$$\begin{aligned} \operatorname{Im} \psi^0(it + \epsilon) &= 2\delta \operatorname{Im} \ln[\cosh(a\epsilon/2 + i(at - b)/2))] \\ &= 2\delta \operatorname{Im} \ln \left[ \frac{e^{a\epsilon/2} + e^{-a\epsilon/2}}{2} \cos \frac{at - b}{2} + i \frac{e^{a\epsilon/2} - e^{-a\epsilon/2}}{2} \sin \frac{at - b}{2} \right]. \end{aligned}$$

If $\cos \frac{at - b}{2} > 0$, $\operatorname{Im} \psi(it \pm 0) = 0$. If $\cos \frac{at - b}{2} < 0$ and $\pm \sin \frac{at - b}{2} > 0$, $\operatorname{Im} \psi(it + 0) = \pm 2\delta\pi$ and $\operatorname{Im} \psi(it - 0) = \mp 2\delta\pi$. Hence, Meixner processes are sSL-processes but not SL-processes. The sSL-measures are $\mathcal{G}_{\pm}(dt) = g_{\pm}(t)dt$, where, for $t > 0$,

$$g_+(t) = 2\delta \sum_{k=0}^{+\infty} \left( \mathbf{1}_{a^{-1}((4k+1)\pi + b, (4k+2)\pi + b)}(t) - \mathbf{1}_{a^{-1}((4k+2)\pi + b, (4k+3)\pi + b)}(t) \right),$$

$$g_-(t) = 2\delta \sum_{-\infty}^{-1} \left( -\mathbf{1}_{a^{-1}((4k+1)\pi + b, (4k+2)\pi + b)}(-t) + \mathbf{1}_{a^{-1}((4k+2)\pi + b, (4k+3)\pi + b)}(-t) \right).$$

## 4. Wiener-Hopf factors

**4.1. General formulas for RLPE.** Let $X, \bar{X}, \underline{X}$ be an RLPE and the supremum and infimum processes of $X$, all starting at 0, and let $T_q, q > 0$, be an exponentially distributed random variable of mean $1/q$, independent of $X$. The Wiener-Hopf factors $\phi_q^+(\xi) = \mathbb{E}[e^{i\xi \bar{X}_{T_q}}]$, $\phi_q^-(\xi) = \mathbb{E}[e^{i\xi \underline{X}_{T_q}}]$ admit analytic continuation to wide regions of $\mathbb{C}^2$, and a number of analytic representations of $\phi_q^{\pm}(\xi)$ can be derived.

Let $q > 0$. In [13, 12] (see also [11, 35]), the following statements are proved.

I. There exist $\sigma_-(q) < 0 < \sigma_+(q)$ such that

$$\tag{4.1} q + \psi(\eta) \notin (-\infty, 0], \quad \operatorname{Im} \eta \in (\sigma_-(q), \sigma_+(q)).$$

II. The Wiener-Hopf factor $\phi_q^+(\xi)$ admits analytic continuation into the half-plane $\operatorname{Im} \xi > \sigma_-(q)$, and can be calculated as follows: for any $\omega_- \in (\sigma_-(q), \operatorname{Im} \xi)$,

$$(4.2) \qquad \phi_q^+(\xi) = \exp\left[\frac{1}{2\pi i} \int_{\operatorname{Im} \eta = \omega_-} \frac{\xi \ln(1 + \psi(\eta)/q)}{\eta(\xi - \eta)} d\eta\right]$$

$$(4.3) \qquad = \exp\left[-\frac{1}{2\pi i} \int_{\operatorname{Im} \eta = \omega_-} \frac{\psi'(\eta)}{q + \psi(\eta)} \cdot \ln \frac{\eta}{\eta - \xi} d\eta\right].$$

III. The Wiener-Hopf factor $\phi_q^-(\xi)$ admits analytic continuation into the half-plane $\operatorname{Im} \xi < \sigma_+(q)$, and can be calculated as follows: for any $\omega_+ \in (\operatorname{Im} \xi, \sigma_+(q))$,

$$(4.4) \qquad \phi_q^-(\xi) = \exp\left[-\frac{1}{2\pi i} \int_{\operatorname{Im} \eta = \omega_+} \frac{\xi \ln(1 + \psi(\eta)/q)}{\eta(\xi - \eta)} d\eta\right]$$

$$(4.5) \qquad = \exp\left[\frac{1}{2\pi i} \int_{\operatorname{Im} \eta = \omega_+} \frac{\psi'(\eta)}{q + \psi(\eta)} \cdot \ln \frac{\eta}{\eta - \xi} d\eta\right].$$

**Remark 4.1.** 1. If $X$ is an sSL process, then $\sigma_-(q) = -\beta_q^+$ if $\beta_q^+$ exists and $\sigma_-(q) = \mu_-$ otherwise, and $\sigma_+(q) = -\beta_q^-$ if $\beta_q^-$ exists and $\sigma_+(q) = \mu_+$ otherwise.

2. To define analytic continuation w.r.t. $q$ into wide regions of the complex plane, it is necessary or advantageous to deform the contours of integration. See [16, 22, 21] for details.

3. If the contours are deformed, then to use the formulas above smaller $\omega_\pm$ may be needed. The key condition is that $q + \psi(\eta) \notin (-\infty, 0]$ or $1 + \psi(\eta)/q \notin (-\infty, 0]$ for all $q$ that appear in the Laplace inversion formula and $\eta \in S_{[\mu_-, \mu_+]}$.

4. In Section 6, a different set of formulas, for $q$ small in absolute value, is used to calculate the asymptotics of $\mathbb{P}[\bar{X}_t < x]$ for $x > 0$ fixed and $t \to \infty$.

5. In [17], formulas for the Wiener-Hopf factors similar to (4.2)-(4.5) are derived for stable Lévy processes. The key differences are 1) the integration is over contours of the form $\Sigma_\gamma = e^{i(\pi-\gamma)} R_+ \cup e^{i\gamma} R_+$ and $\Sigma_{-\gamma} = e^{i(-\pi+\gamma)} R_+ \cup e^{-i\gamma} R_+$, where $\gamma \in (0, \pi/2)$; 2) certain regularization of the integrals is used to ensure the absolute convergence. The generalization of the formulas in [17] for mixtures of stable Lévy processes and processes with exponentially decaying tails of the Lévy densities is straightforward, and then the proofs of the formulas below for the case of large $q$ can be repeated essentially verbatim. However, if the mixture is not SL process, then the problem of location of zeros and study of their analytical properties is highly non-trivial.

Let $\psi$ satisfy

$$(4.6) \qquad \psi(\xi) = c_\pm |\xi|^\nu + O(|\xi|^{\nu-\delta}), \ \xi \to \pm\infty,$$

where $\delta > 0$, $c_\pm = ce^{\pm i\varphi}$, $c > 0$ and $\varphi \in (-\pi/2, \pi/2)$. Define $\nu_\pm = \nu/2 \mp \varphi/\pi$. Straightforward calculations prove the following lemma, which is the basis for the proof of the modified formulas for the Wiener-Hopf factors.

**Lemma 4.1.** *(a) If (4.6) is satisfied,*

$$(4.7) \qquad \psi_{norm}(\xi) := c^{-1}(1 - i\xi)^{-\nu_+}(1 + i\xi)^{-\nu_-}\psi(\xi) = 1 + O(|\xi|^{-\delta}), \ \xi \to \pm\infty.$$

*If the order of the process $\nu < 1$ and $\mu \neq 0$,*

$$(4.8) \qquad \psi_{norm}(\xi) := (1 - i\mu\xi)^{-1}\psi(\xi) = 1 + O(|\xi|^{\nu-1}), \ \xi \to \pm\infty.$$

**Remark 4.2.** Conditions of Lemma 4.1 are satisfied for stable Lévy processes with drift, mixtures of stable Lévy processes with drift and SINH-regular process (the conditions are satisfied and asymptotic formulas hold as $\xi \to \infty$ in appropriate cones in the right and left half-planes). For RLPE (resp., SINH-regular processes), conditions are satisfied and asymptotic formulas hold as $\operatorname{Re} \xi \to \pm\infty$, $\xi$ remaining in appropriate strips (resp., unions of appropriate strips and cones).

Let $X$ be SINH-regular, and let $\sigma > 0$. Then there exists $\omega < 0$ such that for all $q$ satisfying $\operatorname{Re} q = \sigma > 0$, $\phi_q^+(\xi)$ can be represented in several forms depending on the order $\nu$ (and $\mu$, if $\nu < 1$). The proof can be found in [12, Sect.8.2.3]. (The notation in [6, 34] is slightly different.)

(a) If (4.6) is satisfied,

$$(4.9) \qquad \phi_q^\pm(\xi) \;=\; \left( \frac{(q/c)^{1/\nu}}{(q/c)^{1/\nu} \mp i\xi} \right)^{\nu_\pm} \exp[b_q^\pm(\xi) - b_q^\pm(0)], \;\; \pm\operatorname{Im}\xi > -\omega',$$

where, for any $\omega' \in (\omega, 0)$,

$$(4.10) \qquad b_q^\pm(\xi) \;=\; \pm\frac{1}{2\pi i} \int_{\operatorname{Im}\eta = \mp\omega'} \frac{\ln \Phi(q, \eta)}{\eta - \xi} d\eta,$$

$$(4.11) \qquad \Phi(q, \eta) \;=\; \frac{c\left( (q/c)^{1/\nu} - i\eta \right)^{\nu_+} \left( (q/c)^{1/\nu} + i\eta \right)^{\nu_-}}{q + \psi(\eta)}.$$

(b) If $\nu < 1$ and $\mu > 0$,

$$(4.12) \qquad \phi_q^+(\xi) \;=\; \frac{q}{q - i\mu\xi} \exp[b_q^+(\xi) - b_q^+(0)], \;\; \operatorname{Im}\xi > -\omega',$$

$$(4.13) \qquad \phi_q^-(\xi) \;=\; \exp[b_q^-(\xi) - b_q^-(0)], \;\; \operatorname{Im}\xi < \omega',$$

where, for any $\omega' \in (\omega, 0)$,

$$(4.14) \qquad b_q^\pm(\xi) \;=\; \pm\frac{1}{2\pi i} \int_{\operatorname{Im}\eta = \mp\omega'} \frac{\ln \frac{q - i\mu\eta}{q - i\mu\eta + \psi^0(\eta)}}{\eta - \xi} d\eta.$$

(c) If $\nu < 1$ and $\mu < 0$,

$$(4.15) \qquad \phi_q^+(\xi) \;=\; \exp[b_q^+(\xi) - b_q^+(0)], \;\; \operatorname{Im}\xi > -\omega',$$

$$(4.16) \qquad \phi_q^-(\xi) \;=\; \frac{q}{q - i\mu\xi} \exp[b_q^-(\xi) - b_q^-(0)], \;\; \operatorname{Im}\xi < \omega',$$

where $b_q^\pm(\xi)$ are given by (4.14).

In Cases (b) and (c), set $\nu_+ = 1$, $\nu_- = 0$, and $\nu_+ = 0$, $\nu_- = 1$, respectively. Then, as $\xi \to \infty$ in the upper (resp., lower) half-plane and $q$ is fixed, $\phi_q^\pm(\xi)$ enjoys the following asymptotics

$$(4.17) \qquad \phi_q^\pm(\xi) = \phi_{q,\infty}^\pm (1 - i\xi)^{-\nu_\pm} (1 + O(|\xi|^{-\rho'})),$$

where $\rho' > 0$ and $\phi_{q,\infty}^\pm = (q/c)^{\nu_\pm/\nu} \exp[-b_q^\pm(0)]$ in Case (a), and $\phi_{q,\infty}^\pm = (q/\mu)^{\nu_+} \exp[-b_q^\pm(0)]$ in Cases (b),(c). The proof is straightforward; the details can be found in [12, 6, 34].

4.2. **The case of SL and sSL process and** $q > 0$. Let $X$ be an SL-process. Denote by $U_\pm$ the support of the absolute continuous component of the SL-measure $\mathcal{G}_\pm$, and by $\mathcal{Z}_+(q)$ (resp., $\mathcal{Z}_-(q)$) the set of zeros of $q + \psi(\xi)$ in the upper (resp., lower) half-plane. Since $X$ is and SL-process and $q > 0$, $\mathcal{Z}_+ \subset i(0, +\infty)$ and $\mathcal{Z}_- \subset i(-\infty, 0)$.

**Proposition 4.2.** *Let the following conditions hold:* $q > 0$, $X$ *is an SL process, and all zeros (if exist) are simple. Then*

$$(4.18) \qquad \phi_q^+(\xi) = \prod_{-i\beta \in \mathcal{Z}_-(q)} \frac{\beta}{\beta - i\xi} \exp\left[\frac{1}{\pi} \int_{U_+} \ln \frac{z}{z - i\xi} d\arg(q + \psi(-iz - 0))\right],$$

$$(4.19) \qquad \phi_q^-(\xi) = \prod_{-i\beta \in \mathcal{Z}_+(q)} \frac{-\beta}{-\beta + i\xi} \exp\left[\frac{1}{\pi} \int_{U_-} \ln \frac{z}{z + i\xi} d\arg(q + \psi(iz + 0))\right].$$

*Proof.* First, note that if we understand the integrals above as the Riemann-Stieltjes integrals, then we can omit the products in front of the exponential sign and replace $U_\pm$ with $(0, +\infty)$. We prove (4.18); the proof of (4.19) is by symmetry. Since $q + \psi(\eta) \neq 0$ for all $\eta \notin i\mathbb{R}$, we can move the line of integration on the RHS of (4.3) down. If the zero $-i\beta_q^+$ exists, we cross $-i\beta_q^+$ applying the residue theorem, and then deform the contour further, to the banks of the cut $i(-\infty, \mu_-]$. In the remaining integral, we change the variable $\eta = -iz$ and obtain, for $\xi$ in the upper half-plane,

$$(4.20) \qquad \phi_q^+(\xi) = \frac{\beta_q^+}{\beta_q^+ - i\xi} \exp b_q^{+,0}(\xi)$$

(if $-i\beta_q^+$ does not exists, the fraction on the RHS of (4.20) must be omitted), where

$$(4.21) \qquad b_q^{+,0}(\xi) = \frac{1}{2\pi} \int_\infty^{-\mu_-} dz \left[\frac{\psi'(-iz - 0)}{q + \psi(-iz - 0)} - \frac{\psi'(-iz + 0)}{q + \psi(-iz + 0)}\right] \ln \frac{z}{z - i\xi}.$$

Taking into account that $\psi'(-iz \pm 0) = id(q + \psi(-iz \pm 0))/dz$ and $\arg(q + \psi(-iz + 0)) = -\arg(q + \psi(-iz - 0))$, we obtain

$$(4.22) \qquad b_q^{+,0}(\xi) = \frac{1}{\pi} \int_{-\mu_-}^\infty \ln \frac{z}{z - i\xi} d\arg(q + \psi(-iz - 0)).$$

$\square$

**Remark 4.3.** 1. In the case of SL-processes, all zeros are on $i\mathbb{R}$. In the case of sSL-processes, there can be (probably, always will be, but we do not have a proof) zeros outside $i\mathbb{R}$. For instance, if $X$ is a "driftless" Meixner process, then there are two infinite series of zeros (depending on $q$); the calculation of the roots is straightforward. The proof of (4.18) and (4.19) remains valid if we assume that the sets of zeros do not have accumulation points, which is the case for the "driftless" Meixner process.

2. To study the survival probability problem, we need formulas for $\phi_q^\pm(\xi)$ when $q \in \mathbb{C}$ is small and asymptotic formulas as $q \to 0$. The formulas and properties depend on the sign of the first instantaneous moment $\mu_1 := i\psi'(0)$. We will derive these formulas in Sect. 6 for Lévy processes with exponentially decaying tails. In Sections 6.3 and 6.5, we refine the formulas for the asymptotic coefficients for SL-processes in terms of zeros of $\psi$ and the integral w.r.t. to the absolute continuous components of the SL-measure.

3. Similar formulas for the lower tail probability problem are possible to derive but only assuming that zeros of $q + \psi(\xi)$ for $q$ with a sufficiently large real part are known. This is the case for 1) a "driftless" Meixner model, 2) HEJD and $\beta$-models; 3) mixtures of BM and KoBoL or NTS or VGP. In the case of a meromorphic model, additional conditions are need. Hence, a general results are impossible to formulate in even moderately succint fashion although the idea of the proof is the same as in the case of real and small $q$.

## 5. Lower tail probability problem

5.1. **General outline.** For regular SINH processes, the lower tail probability problem is solved in [6, 34] in the option pricing framework: if the riskless rate is 0, then the price of the no-touch barrier option with an upper barrier is the survival probability. The results are obtained for RLPEs satisfying additional unnecessary stringent conditions. In fact, the proof uses only the properties of SINH-regular processes. We reformulate the results of [6, 34] for the leading term of the asymptotics of $\mathbb{P}[\bar{X}_t < x]$ for $t$ fixed and $x \downarrow 0$. Note that in [6], the second term of asymptotics is calculated as well, and the reformulation of this result for the lower tail probability problem is possible. The asymptotic formulas with the leading term only are of the form

$$\mathbb{P}[\bar{X}_t < x] \sim \kappa(t)x^{\nu_+}, \ x \downarrow 0, \tag{5.1}$$

where $\nu_+ \in [0, 1]$ is the opposite to the rate of decay of the positive Wiener-Hopf factor $\phi_q^+(\xi)$ as $\xi \to \infty$ in the upper half plane, and $\kappa(t) > 0$ is calculated in terms of the coefficient $\phi_{q,\infty}^+$ in the asymptotic formula (4.17) for $\phi_q^+(\xi)$. The exponent $\nu_+ = 0$ appears when, for $q > 0$, the distribution of $\bar{X}_{T_q}$ has an atom at 0. We reformulate some of the results obtained in [6, 34] for barrier options. The starting point is the standard formula

$$\mathbb{P}[\bar{X}_t < x] = \frac{1}{2\pi i} \int_{\mathrm{Re}\, q = \sigma} dq \, \frac{e^{qt}}{q} \frac{1}{2\pi} \int_{\mathrm{Im}\, \xi = \omega} d\xi \, e^{-ix\xi} \frac{\phi_q^+(\xi)}{i\xi}, \tag{5.2}$$

where $\sigma > 0$ and $\omega < 0$ are sufficiently large and small, respectively, so that $\mathrm{Re}(q + \psi(\xi))$ is bounded away from 0 for all $q \in \{q \mid \mathrm{Re}\, q \geq \sigma\}$ and $\xi \in S_{[\omega,0]}$. In [6, Sect. 1.1] and Theorem 6.3 in [34], equations similar to (5.2) are used to derive the leading (and, under additional conditions, second) term of the asymptotics as $x \downarrow 0$. The results are formulated for processes with the positive killing rate $r$ (riskless rate) but the proof is valid for $r = 0$ as well. Next, the results are formulated for the expectation of $\mathbf{1}_{(x,+\infty)}(\underline{X}_t)G(X_t)$, where $x < 0$, and $G$ is a sufficiently regular function. In our case, $G(x) = 1$ for all $x$, hence, the function $\tilde{G}(q, 0+)$ in part (a) of Theorem 6.3 in [34] equals 1. Applying the standard symmetry considerations, we reformulate the results as follows.

**Theorem 5.1.** *Let $X$ be SINH-regular of order $\nu \in (0, 2]$ or of order $0+$ with negative drift. Then there exists $\sigma > 0$ such that*

(a) $q \mapsto \phi_{q,\infty}^+$ *is analytic in the half-plane* $\mathrm{Re}\, q \geq \sigma$;

(b) *there exist* $C, s > 0$ *and integer* $k \geq 0$ *such that*

$$|\partial_q^k(q^{-1}\phi_{q,\infty}^+)| \leq C|q|^{-1-s}, \quad \mathrm{Re}\, q \geq \sigma, \tag{5.3}$$

*hence, the following integral absolutely converges:*

$$(5.4) \qquad \kappa_k(t) = \frac{1}{2\pi i \Gamma(1+\nu_+)(-t)^k} \int_{\mathrm{Re}\, q=\sigma} e^{qT} \partial_q^k(q^{-1}\phi_{q,\infty}^+) dq;$$

*(c) there exist an integer $k \geq 0$ and $s > 0$ such that as $x \downarrow 0$,*

$$(5.5) \qquad \mathbb{P}[\bar{X}_t < x] = \kappa_k(t)x^{\nu_+} + O(x^{\nu_++s}).$$

**Remark 5.1.** 1. If $\nu_+ > 0$, then $k = 0$ is admissible.

2. If we understand the integral in principal value sense, $k = 0$ is admissible in all cases.

3. If $\nu_+ = 0$, $\lim_{x\downarrow 0} \mathbb{P}[\bar{X}_t < x] > 0$.

4. Using the standard symmetry considerations, we can rewrite (5.4) in the form

$$\kappa_0(t) = \frac{1}{\pi \Gamma(1+\nu_+)} \lim_{A \to +\infty} \mathrm{Im} \int_\sigma^{\sigma+iA} e^{qT} q^{-1} \phi_{q,\infty}^+ dq.$$

5. If, for $q$ with a large Re $q$, the zeros of $q + \psi(\xi)$ can be calculated, then $\phi_{q,\infty}^+$ can be calculated in terms of the zeros and integrals over $U_+$. See Remark 4.3.4.

## 6. Asymptotics of the survival probability as $t \to \infty$

6.1. **General outline.** We prove that depending on the sign of the first instantaneous moment $\mu_1 := i\psi'(0)$, 3 possibilities exist. In cases $\mu_1 < 0, \mu_1 > 0$ and $\mu_1 = 0$, respectively,

$$(6.1) \qquad \mathbb{P}[\bar{X}_t < x] = 1 - p_\infty(x) + O(t^{-1}),$$
$$(6.2) \qquad \mathbb{P}[\bar{X}_t < x] = O(t^{-1}),$$
$$(6.3) \qquad \mathbb{P}[\bar{X}_t < x] = p_{\infty,0}(x)t^{-1/2} + O(t^{-1}),$$

where $p_\infty(x) \in (0,1)$ and $p_{\infty,0}(x) > 0$. The formulas for the asymptotic coefficients are refined for SL-processes. We derive analytic expressions for $p_\infty(x)$ and $p_{\infty,0}(x)$ in terms of $\mu_1$ and $\mu_2 := \psi''(0)$, respectively, zeros of $\psi$, the restriction of $\phi_q^-(\xi)$ on $i(-\infty, 0)$, and densities of the absolute continuous components of the SL-measures. In all cases, we formulate the results assuming that the positive and negative jump components are non-trivial, and zeros of $\psi$ are simple. The reader can easily simplify the results (and the proof) for the case when one jump component is absent (or both) and some of zeros are of multiplicity greater than 1. The resulting formulas and proofs can be modified for sSL-processes provided the location of zeros of $\psi$ is known.

Assuming that the asymptotic formulas for the left tail probability in the case of mixtures of stable Lévy processes are derived, one may try, as in the case of stable Lévy processes, to use the rescaling and reduce the survival probability problem for mixtures of stable Lévy processes and SINH-regular processes to the lower tail probability problem. However, in the case of mixtures of stable Lévy processes and regular SINH-processes, the modification of the technique of [6, 34, 15, 17] is not straightforward although possible.[3]

---

[3]In [6, 34], the proof is based on the localization in the dual $(q, \xi)$-space: bounds and asymptotics of the integrand in the Laplace-Fourier inversion formulas are different in appropriately chosen regions in the $(q, \xi)$-space. If $\psi$ depends on $t$, it is necessary to localize in the $(t, q, \xi)$-space (microlocalize), and a modification of the microlocalization technique developed in [31] (see also [36]) to study boundary problems for degenerate elliptic operators can be applied.

6.2. **General formulas: first steps.** Let $X$ be an RLPE of order $\nu \in (0, 2]$. Then, as special cases of general theorems for barrier options in [12, 11]), we have the following two formulas. For any $\sigma_0 > 0$, there exists $\omega_- < 0$ such that 1) $q + \psi(\xi) \notin (-\infty, 0]$ for all $q$ on the line $\{q \mid \mathrm{Re}\, q = \sigma_0\}$ and all $\xi \in S_{(\omega', 0]}$, and 2) for any $x > 0$ and $t > 0$,

$$(6.4) \qquad \mathbb{P}[\bar{X}_t \geq x] = \frac{1}{2\pi i} \int_{\mathrm{Re}\, q = \sigma_0} dq\, \frac{e^{qt}}{q}\, \frac{1}{2\pi} \int_{\mathrm{Im}\, \xi = \omega_-} d\xi\, e^{-ix\xi} \frac{\phi_q^+(\xi)}{i\xi},$$

$$(6.5) \qquad \mathbb{P}[\bar{X}_t < x] = \frac{1}{2\pi i} \int_{\mathrm{Re}\, q = \sigma_0} dq\, \frac{e^{qt}}{q}\, \left(1 - \frac{1}{2\pi} \int_{\mathrm{Im}\, \xi = \omega_-} d\xi\, e^{-ix\xi} \frac{\phi_q^+(\xi)}{i\xi}\right).$$

In [7], the same formulas are proved for any Lévy process with the characteristic exponent analytic in a strip around the real axis. Unless otherwise stated, the formulas that we derive below and proofs are valid for any process satisfying this condition. Additional assumptions such as $X$ is an RLPE or SL process are imposed whenever necessary.

Using the Wiener-Hopf factorization identity, we rewrite (6.5) as

$$(6.6) \qquad \mathbb{P}[\bar{X}_t < x] = \frac{1}{2\pi i} \int_{\mathrm{Re}\, q = \sigma_0} dq\, \frac{e^{qt}}{q}\, \left(1 - \frac{1}{2\pi} \int_{\mathrm{Im}\, \xi = \omega_-} \frac{q e^{-ix\xi} d\xi}{i\xi \phi_q^-(\xi)(q + \psi(\xi))}\right).$$

We use (6.4) or (6.6) if $\mu_1 < 0$ or $\mu_1 > 0$, respectively. If $\mu_1 = 0$, we impose an additional condition $\nu_+ > 0$. Then $\lim_{x \downarrow 0} \mathbb{P}[X_t \geq x] = 1$ and the integral on the RHS of (6.4) converges absolutely and uniformly in $x \in (0, 1]$. Therefore, we can pass to the limit $x \downarrow 0$ in (6.4) and, using the Wiener-Hopf factorization, obtain

$$(6.7) \qquad 1 = \mathbb{P}[\bar{X}_t \geq 0] = \frac{1}{2\pi i} \int_{\mathrm{Re}\, q = \sigma_0} dq\, \frac{e^{qt}}{q}\, \frac{1}{2\pi} \int_{\mathrm{Im}\, \xi = \omega_-} \frac{q\, d\xi}{i\xi \phi_q^-(\xi)(q + \psi(\xi))}.$$

Substituting (6.4) and (6.7) into $\mathbb{P}[\bar{X}_t < x] = 1 - \mathbb{P}[\bar{X}_t \geq 0]$, we obtain

$$(6.8) \qquad \mathbb{P}[\bar{X}_t < x] = \frac{1}{2\pi i} \int_{\mathrm{Re}\, q = \sigma_0} dq\, \frac{e^{qt}}{q}\, \frac{1}{2\pi} \int_{\mathrm{Im}\, \xi = \omega_-} \frac{q(1 - e^{-ix\xi}) d\xi}{i\xi \phi_q^-(\xi)(q + \psi(\xi))}.$$

For a large $t$, we choose $\sigma_0$ in the form $\sigma_0 = \sigma/t$, where $\sigma > 0$ is independent of $t$. In view of Lemma 6.2, if $\mu_1 > 0$, $\sigma > 0$ must be sufficiently small, and if $\mu_1 = 0$, we can choose any $\sigma > 0$ and the constructions below are valid for sufficiently large $t$. Finally, if $\mu_1 < 0$, then we may fix a sufficiently small $\sigma_0 > 0$ independent of $t$ but we need to use $\sigma_0 = \sigma/t$ so that $|e^{qt}|$ is uniformly bounded on the line of integration. Considering separately cases 1) $\nu \in [1, 2]$; 2) $\nu \in (0, 1)$, $|q - i\mu\xi| \leq |\xi|^\nu$; and 3) $\nu \in (0, 1)$, $|q - i\mu\xi| > |\xi|^\nu$; we derive the bound

$$(6.9) \qquad |\partial_q^k (q + \psi(\xi))^{-1}| \leq C(|q| + |\xi|)^{-(1+k)\min\{1,\nu\}}, k = 0, 1, \ldots.$$

Integrating by parts and using (6.9), we obtain that the integrals (6.4), (6.6) and (6.8) over $\{q \mid \mathrm{Re}\, q = \sigma/t, |\mathrm{Im}\, q| \geq \epsilon\}$ are $O(t^{-1})$. It remains to calculate the asymptotics of the integrals over $\mathcal{L}_{\sigma, \epsilon; t} = \{q \mid \mathrm{Re}\, q = \sigma/t, |\mathrm{Im}\, q| \leq \epsilon\}$. Denote these integrals $V_{\sigma, \epsilon}^+(t, x)$, $V_{1, \sigma, \epsilon}^-(t, x)$ and $V_{2, \sigma, \epsilon}^-(t, x)$, respectively. On the strength of Lemma 6.2, we can find $t_0$ and $\omega'_-(< \omega_-)$ such that for each $q \in \mathcal{L}_{\sigma, \epsilon; t}$ and $t \geq t_0$, 1) $q + \psi(\eta) \neq 0$ for all $\eta \in S_{[\omega'_-, 0]}$ in the case $\mu_1 < 0$, and 2) in the case $\mu_1 \geq 0$, the equation $q + \psi(\xi) = 0$ has the only root $-i\beta_q^+$ in $S_{(\omega'_-, \omega_-]}$, and the root is simple.

6.3. **The case** $\mu_1 < 0$. In this case, we can choose $\epsilon > 0$ and $\omega' < \omega < 0$ so that, for each $q \in B(0, \epsilon)$, the root $-i\beta_q^+$ either does not exist or lies below the line $\{\mathrm{Im}\,\eta = \omega'\}$ and $q + \phi(\eta) \neq 0$ for all $\eta \in S_{[\omega', 0]}$. Hence, we can define $\phi_q^+(\xi)$ for $\xi \in S_{(\omega'_+, +\infty)}$ above this line by (4.3) with $\omega'_-$ in place of $\omega_-$. It follows that for all $\xi \in S_{(\omega'_+, +\infty)}$, $\phi_q^+(\xi)$ is analytic function of $q \in B(0, \epsilon)$. Therefore, deforming the contour $\mathcal{L}_{\sigma,\epsilon;t}$ on the RHS of (6.4) into $\mathcal{L}'_{\sigma,\epsilon;t} = \{q \mid |q| = \epsilon, \mathrm{Re}\,q \leq \sigma/t\}$, crossing the simple pole at $q = 0$ and applying the residue theorem, we obtain

$$(6.10) \quad V_{\sigma,\epsilon}^+(t, x) = \frac{1}{2\pi}\int_{\mathrm{Im}\,\xi = \omega_-} d\xi\, e^{-ix\xi}\frac{\phi_0^+(\xi)}{i\xi} + \frac{1}{2\pi i}\int_{\mathcal{L}'_{\sigma,\epsilon;t}} dq\, \frac{e^{qt}}{q}\frac{1}{2\pi}\int_{\mathrm{Im}\,\xi = \omega_-} d\xi\, e^{-ix\xi}\frac{\phi_q^+(\xi)}{i\xi}.$$

Integrating by parts in the outer integral in the second term on the RHS of (6.10), we conclude that the second term is $O(t^{-1})$ and obtain

**Theorem 6.1.** *Let $\mu_1 < 0$ and $x > 0$. Then (6.1) holds with*

$$(6.11) \qquad\qquad p_\infty(x) = \frac{1}{2\pi}\int_{\mathrm{Im}\,\xi = \omega'_-} d\xi\, e^{-ix\xi}\frac{\phi_0^+(\xi)}{i\xi},$$

*where $\omega'_- < 0$ satisfying $\psi(\xi) \neq 0$ for $\xi \in S_{[\omega'_-, 0)}$ is arbitrary.*

To formulate a refined version of (6.11) for SL processes, we need a modification of the formula for $\phi_q^-(\xi)$ in the case $\mu_1 \leq 0$, for $\xi$ in the lower half-plane, where $q$ is in a ball $B(0, \epsilon)$ of a small radius $\epsilon > 0$, centered at 0. The starting point is the following lemma.

**Lemma 6.2.** *Let $\psi$ be analytic in a strip around the real axis, and let $\mathrm{Re}\,\psi(\xi) \to \infty$ as $\xi \to \infty$ in the strip. Then there exist $\epsilon, \epsilon_1 > 0$ and $\omega' > 0$ such that*

*1) for all $q \in B(0, \epsilon)$, the equation $q + \psi(\xi) = 0$ has the following zeros in the ball $B(0, \epsilon_1)$:*

(a) *if $\mu_1 > 0$, one zero (of multiplicity 1) $-i\beta_q^+ \in B(0, \epsilon_1)$; the function $q \mapsto \beta_q^+$ is analytic in $B(0, \epsilon_1)$, and $\lim_{q \to 0} \beta_q^+/q = 1/\mu_1 > 0$;*

(b) *if $\mu_1 < 0$, one zero (of multiplicity 1) $-i\beta_q^- \in B(0, \epsilon_1)$; the function $q \mapsto \beta_q^-$ is analytic in $B(0, \epsilon_1)$, and $\lim_{q \to 0} \beta_q^-/q = 1/\mu_1 < 0$;*

(c) *if $\mu_1 = 0$ and $q \notin (-\omega, 0]$, two zeros $\pm i\beta_q^\pm$ (each of multiplicity one). Functions $q \mapsto \beta_q^\pm$ are analytic function of $z := \sqrt{q}$ on $\{z \mid |z| \leq \sqrt{\epsilon}, \mathrm{Re}\,z \geq 0\}$ (meaning analytic in the interior and continuous at the boundary) and $\beta_q^\pm/\sqrt{q} \to \pm\sqrt{2/\mu_2}$, as $q \to 0$ remaining in $\mathbb{C} \setminus (-\infty, 0]$;*

*2) there exists $c(\epsilon_1, \omega') > 0$ such that $|q + \psi(\xi)| \geq c(\epsilon_1, \omega')$ for all $S_{(-\omega', \omega')} \setminus B(0, \epsilon_1)$.*

*Proof.* We use the asymptotic expansion $\psi(\xi) = \frac{\mu_2}{2}\xi^2 - i\mu_1\xi + \cdots$ as $\xi \to 0$. (a) Since $\psi'(0) = -i\mu_1 \neq 0$, $\psi'(\xi) \neq 0$ in a neighborhood of 0, hence, the implicit function theorem is applicable, and $\beta_q^+$ is an analytic function in the neighborhood of 0 satisfying $\lim_{q \to 0} \beta_q^+/q = 1/\mu_1 > 0$. (b) the same proof. (c) Since $\partial\psi(\xi)/d(\xi^2) = \mu_2/2 > 0$, the same argument proves that $\beta_q^\pm$ are analytic functions of $z = \sqrt{q}$ in $\{z \mid |z| \leq \sqrt{\epsilon}, \mathrm{Re}\,z > 0\}$ and continuous at the boundary, and $\beta_q^\pm/\sqrt{q} \to \pm\sqrt{2/\mu_2}$, as $q \to 0$ remaining in $\mathbb{C} \setminus (-\infty, 0]$. 2) holds since $\mathrm{Re}\,\psi(\xi) > 0$ for all $\xi \in \mathbb{R} \setminus 0$ and $|\psi(\xi)| \to \infty$ as $\xi \to \infty$ in the strip. $\qquad\square$

We use the following modification of the formula for $\phi_q^-(\xi)$, which is similar to several modifications of the formulas for the Wiener-Hopf factors in [35] and can be easily verified:

$$(6.12) \qquad \phi_q^-(\xi) = \frac{-\beta_q^-}{-\beta_q^- + i\xi} \phi_q^{-,0}(\xi),$$

where

$$(6.13) \qquad \phi_q^{-,0}(\xi) = \exp\left[-\frac{1}{2\pi i}\int_{\operatorname{Im}\eta=\omega_+} \frac{\xi\ln\Phi(q,\eta)d\eta}{\eta(\eta-\xi)}\right],$$

$$(6.14) \qquad \Phi(q,\eta) = \frac{q(-\beta_q^- + i\eta)}{-\beta_q^-(q+\psi(\eta))},$$

and $\omega_+ > 0$ is sufficiently small. By continuity, define $\Phi(0,\eta) = -i\mu_1\eta/\psi(\eta)$. The following lemma is immediate from Lemma 6.2.

**Lemma 6.3.** *Let $\mu_1 \leq 0$. Then there exists $\epsilon > 0$ and $\omega_+ > 0$ such that for all $q \in B(0,\epsilon)$ and $\eta \in S_{(0,\omega_+]}$,*

*(a) $\Phi(q,\eta) \neq 0$, hence, $\ln\Phi(q,\eta)$ is a well-defined function on an appropriate Riemann surface satisfying the requirement $\ln 1 = 0$;*

*(b) $\phi_q^{-,0}(\xi)$ and $1/\phi_q^{-,0}(\xi)$ and their derivatives w.r.t. $q$ are bounded analytic functions of $q \in B(0,\epsilon)$, if $\mu_1 < 0$, and of $z := \sqrt{q} \in \{z \mid \operatorname{Re}z \geq 0, |z| < \sqrt{\epsilon}\}$, if $\mu_1 = 0$.*

We also need

**Proposition 6.4.** *Let $X$ be an SL process, and let all zeros of $\psi$ on $i(0,+\infty)$ (if exist) be simple. Then*

$$(6.15) \qquad \phi_0^{-,0}(\xi) = \prod_{-i\beta\in\mathcal{Z}_+(0),\beta\neq\beta_q^-} \frac{-\beta}{-\beta+i\xi}\exp\left[\frac{1}{\pi}\int_{U_-}\ln\frac{z}{z+i\xi}d\arg\psi(iz+0)\right].$$

*Proof.* The proof is the straightforward modification of the proof of Proposition 4.2. $\qquad\square$

**Theorem 6.5.** *Let the following conditions hold:*

*(a) $X$ is an SL-process and $\mu_1 < 0$;*

*(b) $\int_{U^+} |\psi(-iw-0)|^{-2}\mathcal{G}_+(dw) < \infty$;*

*(c) the zeros of $\psi$ on $i(-\infty,0)$ are simple.*

*Then, for $x > 0$, (6.1) holds with*

$$(6.16) \qquad p_\infty(x) = -\frac{\mu_1}{\pi}\int_{U_+}\frac{e^{-xw}\mathcal{G}_+(dw)}{\phi^{0,-}(-iw)|\psi(-iw-0)|^2} - \mu_1\sum_{-i\beta\in\mathcal{Z}_-(0)\setminus 0}\frac{e^{-x\beta}}{\phi_0^{-,0}(-i\beta)i\psi'(-i\beta)}.$$

*Proof.* Using the Wiener-Hopf factorization formula and (6.12), we obtain

$$\frac{\phi_0^+(\xi)}{i\xi} = \lim_{q\to 0}\frac{q(-\beta_q^+ + i\xi)}{i\xi(-\beta_q^+)\phi_q^{-,0}(\xi)(q+\psi(\xi))} = \frac{-\mu_1}{\phi_0^{-,0}(\xi)\psi(\xi)}.$$

Substituting into (6.11), deforming the contour to the cut $i(-\infty,\mu_-]$, taking into account the contribution of simple poles (zeros of $\psi$), and changing the variable $\xi = -iw \pm 0$ on the left and right banks of the cuts along the connected components of $-U_+$, we obtain

$$
\begin{aligned}
p_\infty(x) \ = \ &-\mu_1 \sum_{-i\beta \in \mathcal{Z}_-(0)\setminus 0} \frac{e^{-x\beta}}{\phi_0^{-,0}(-i\beta)i\psi'(-i\beta)} \\
&+ \frac{-\mu_1}{2\pi} \int_{U_+} dw \, \frac{e^{-wx}}{\phi_0^{-,0}(-iw)} i \left[ \frac{1}{\psi(-iw-0)} - \frac{1}{\psi(-iw+0)} \right].
\end{aligned}
$$

It remains to recall that $-i[\psi(-iw-0) - \psi(-iw+0)]dw = \mathcal{G}_+(dw)$ (see Theorem 3.16, b). $\quad\square$

6.4. **The case $\mu_1 > 0$.** We can move the line of integration in the inner integral on the RHS of (6.6) down to a line $\{\xi \mid \operatorname{Im}\xi = \omega'_-\}$, where $i\omega'_-$ is in the strip of analyticity but below $-i\beta_q^+$ for all $q \in B(0,\epsilon)$. On crossing the simple pole at $\xi = -i\beta_q^+$, we apply the residue theorem:

$$
\begin{aligned}
(6.17) \qquad V_{1,\sigma,\epsilon}^-(t,x) \ = \ & \frac{1}{2\pi i} \int_{\mathcal{L}_{\sigma,\epsilon;t}} dq \, \frac{e^{qt}}{q} \left( 1 - \frac{e^{-x\beta_q^+} q}{\beta_q^+ \phi_q^-(-i\beta_q^+) i\psi'(-i\beta_q^+)} \right) \\
& - \frac{1}{2\pi i} \int_{\mathcal{L}_{\sigma,\epsilon;t}} dq \, e^{qt} \frac{1}{2\pi} \int_{\operatorname{Im}\xi=\omega'_-} d\xi \, \frac{e^{-ix\xi}}{i\xi\phi_q^-(\xi)(q+\psi(\xi))}.
\end{aligned}
$$

If $\mu_1 > 0$, $\beta_q^+ \to q/\mu_1$ as $q \to 0$, and $i\psi'(\xi) \sim \mu_1$, $\phi_q^-(\xi) \to 1$ as $\xi \to 0$. Hence, the function

$$
q \mapsto q^{-1} \left( 1 - \frac{e^{-x\beta_q^+} q}{\beta_q^+ \phi_q^-(-i\beta_q^+) i\psi'(-i\beta_q^+)} \right)
$$

is an analytic function of $q \in B(0,\epsilon)$, vanishing at 0. Therefore, we may deform the contour of integration in the first integral on the RHS of (6.17) into the contour $\mathcal{L}'_{\sigma,\epsilon;t}$, and, integrating by parts, prove that the first term on the RHS of (6.17) is $O(t^{-1})$. The interior integral in the double integral on the RHS of (6.17) is an analytic function of $q \in B(0,\epsilon)$, hence, we may use the same argument to conclude that the second term on the RHS of (6.17) is $O(t^{-1})$ as well. We have proved

**Theorem 6.6.** Let $\mu_1 > 0$ and $x > 0$. Then (6.2) holds.

6.5. **The case $\mu_1 = 0$.** In the formula for $V_{2,\sigma,\epsilon;t}^-$, which is (6.8) with the integration over $\mathcal{L}_{\sigma,\epsilon;t}$ instead of $\{q \mid \operatorname{Re} q = \sigma/t\}$, we move the line of integration down to a line $\{\xi \mid \operatorname{Im}\xi = \omega'_-\}$, where $\omega'_- > \mu_-$ and $i\omega'_-$ is below all $-i\beta_q^+, q \in B(0,\epsilon)$. Crossing the simple pole at $\xi = -i\beta_q^+$, we obtain

$$
\begin{aligned}
(6.18) \qquad V_{2,\sigma,\epsilon}^-(t,x) \ = \ & \frac{1}{2\pi i} \int_{\mathcal{L}_{\sigma,\epsilon;t}} dq \, \frac{e^{qt}}{q} \frac{(1-e^{-x\beta_q^+})q(-\beta_q^- + \beta_q^+)}{\beta_q^+(-\beta_q^-)\phi_q^{-,0}(-i\beta_q^+)i\psi'(-i\beta_q^+)} \\
(6.19) \qquad & + \frac{1}{2\pi i} \int_{\mathcal{L}_{\sigma,\epsilon;t}} dq \, \frac{e^{qt}}{-\beta_q^-} \frac{1}{2\pi} \int_{\operatorname{Im}\xi=\omega'_-} d\xi \, \frac{(1-e^{-ix\xi})(-\beta_q^- + i\xi)}{i\xi\phi_q^{-,0}(\xi)(q+\psi(\xi))}.
\end{aligned}
$$

As $((-\infty, 0] \not\ni)q \to 0$,

$$
\begin{aligned}
\beta_q^+ &= (2q/\mu_2)^{1/2} + c_1^+ q + c_{3/2}^+ q^{3/2} \dots, \\
\frac{-\beta_q^- + \beta_q^+}{\beta_q^+(-\beta_q^+)} &= (2\mu_2/q)^{1/2} + d_0 + d_{1/2} q^{1/2} + \cdots, \\
\phi_q^{-,0}(-i\beta_q^+) &= 1 + u_{1/2} q^{1/2} + \cdots, \\
i\psi(-i\phi_q^+) &= (2q/\mu_2)^{1/2} + v_1 q + \cdots,
\end{aligned}
$$

where $c_1^+ \dots, d_0 \dots, u_{1/2}, \dots, v_1 \dots$, are constants, and the series in powers of $z := q^{1/2}$ converge in a small vicinity of $z = 0$ in the closed right half plane. It follows that the expression in the brackets on the RHS of (6.18) is of the form $x((2q/\mu_2)^{1/2} + (w_0 + w_1 x)q + f(t, x, q)$, where $w_0, w_1 \in \mathbb{C}$, and $f(t, x, q) = O(q^{3/2})$ as $((-\infty, 0] \ni)z \to 0$, uniformly in $t \to \infty$ and $x \in (0, 1]$. Therefore, the first term on the RHS of (6.18) equals

$$
x\left(\frac{2}{\mu_2}\right)^{1/2} \frac{1}{2\pi i} \int_{\mathcal{L}_{\sigma,\epsilon;t}} dq\, e^{qt} q^{-1/2} + (w_1 x + w_2 x^2) \frac{1}{2\pi i} \int_{\mathcal{L}_{\sigma,\epsilon;t}} dq\, e^{qt} + \frac{1}{2\pi i} \int_{\mathcal{L}_{\sigma,\epsilon;t}} dq\, e^{qt} f(t, x, q).
$$

Integrating by parts in the second and third integrals, we obtain that both are $O(t^{-1})$. Deforming the contour of integration, we reduce the first integral above to the sum of two integrals over two parts of the circle: $\{q \mid |q| = \epsilon, \operatorname{Re} q \le \sigma/t, q \ne -\epsilon\}$ (integrating by parts, we prove that both integrals are $O(t^{-1})$), and the integral over the cut:

$$
\begin{aligned}
\frac{1}{2\pi i} \int_{\mathcal{L}_{\sigma,\epsilon;t}} dq\, e^{qt} q^{-1/2} &= \frac{1}{2\pi i} \int_{-\epsilon}^0 e^{zt} \left((iz - 0)^{-1/2} - (iz + 0)^{-1/2}\right) dz + O(t^{-1}) \\
&= \frac{1}{\pi} \int_{-\infty}^0 e^{zt} |z|^{-1/2} \frac{e^{i\pi/2} - e^{-i\pi/2}}{2i} dz + O(t^{-1}) \\
&= t^{-1/2} \frac{\Gamma(1/2)}{\pi} + O(t^{-1}) = t^{-1/2} \pi^{-1/2} + O(t^{-1}).
\end{aligned}
$$

Thus, the term on the RHS of (6.18) admits the asymptotics $xt^{-1/2}(2/(\mu_2\pi))^{1/2} + O(t^{-1})$.

Since $-\beta_q^-$ admits the asymptotic expansion similar to the one for $\beta_q^+$: $-\beta_q^- = (2q/\mu_2)^{1/2} + c_1^- q + \cdots$, $\psi_q^{-,0}(\xi)$ admits the asymptotic expansion of the form $\psi_q^{-,0}(\xi) = \psi_0^{-,0}(\xi) + f_{1/2}(q)q^{1/2} + \cdots$, and $|\psi(\xi)|$ and $|\psi_0^{-,0}(\xi)|$ are bounded away from 0 for $\xi$ on the line of integration, essentially the same argument as above proves that the integral on the RHS of (6.19) admits the asymptotics

$$
t^{-1/2}(2/(\mu_2\pi))^{1/2} \frac{1}{2\pi} \int_{\operatorname{Im}\xi=\omega_-'} \frac{(1 - e^{-ix\xi})d\xi}{\phi_0^{-,0}(\xi)\psi(\xi)} + O(t^{-1}).
$$

We have proved

**Theorem 6.7.** *Let $\mu_1 = 0$ and $x > 0$. Then (6.3) holds with*

$$
\tag{6.20} p_{\infty,0}(x) = \left(\frac{2}{\mu_2\pi}\right)^{1/2} \left(x + \frac{1}{2\pi} \int_{\operatorname{Im}\xi=\omega_-'} \frac{(1 - e^{-ix\xi})d\xi}{\phi_0^{-,0}(\xi)\psi(\xi)}\right),
$$

*where $\omega_-' < 0$ satisfying $\psi(\xi) \ne 0$ for $\xi \in S_{[\omega_-', 0)}$ is arbitrary.*

The following theorem is deduced from Theorem 6.7 similarly to the proof of Theorem 6.5.

**Theorem 6.8.** *Let the following conditions hold:*

*(a) $X$ is an SL-process and $\mu_1 = 0$;*

*(b) $\int_{U_+} |\psi(-iw - 0)|^{-2}\mathcal{G}_+(dw) < \infty$;*

*(c) the zeros of $\psi$ on $i(-\infty, 0)$ are simple.*

*Then, for $x > 0$, (6.1) holds with*

$$p_\infty(x) = \left(\frac{2}{\mu_2 \pi}\right)^{1/2} \left(x + \sum_{-i\beta \in \mathcal{Z}_-(0)\setminus 0} \frac{1 - e^{-x\beta}}{\phi_0^{-,0}(-i\beta)i\psi(-i\beta)} + \frac{1}{\pi}\int_{U_+} \frac{(1 - e^{-xw})\mathcal{G}_+(dw)}{\phi^{0,-}(-iw)|\psi(-iw - 0)|^2}\right).$$

## 7. Conclusion

In the paper, we derived the leading terms of asymptotics of the lower tail probability $\mathbb{P}[\bar{X}_t < x]$ as $x \to 0$, $t$ being ixed, and the asymptotics of the survival function $\mathbb{P}[\bar{X}_t < x]$ as $t \to \infty$, $x > 0$ being fixed. The first problem is solved for wide classes of SINH-regular processes (the second term of asymptotics can be calculated as well). The asymptotic formulas for the second problem are derived for Lévy processes with exponentially decaying tails of the Lévy density. Then, for Stieltjes-Lévy processes (SL-processes), the asymptotic coefficients are expressed in terms of the first two instantaneous moments, zeros of the characteristic exponent $\psi$ and integrals over the supports of absolutely continuous components of Stieltjes-Lévy measures $\mathcal{G}_\pm$ (SL-measures), which are the main ingredients of construction of the process. The formulas can be generalized to the case of signed SL processes (sSL-processes). The results are fairly straightforward consequences of the general formulas for $\mathbb{P}[\bar{X}_t < x]$ and properties of SINH-regular and (s)SL processes. The definitions of SINH-regular and (s)SL-regular processes either explicitly state or allow one to deduce from the definition the key properties of the characteristic exponent which allowed us to develop efficient numerical methods for evaluation of expectations of functions of $X_t$, $(X_t, \bar{X}_t)$, $(X_t, \underline{X}_t)$, joint distributions of $X_t$ and its extremum; joint distribution of $X_t$ and its two extrema; joint pdf of a Lévy process, its extremum, and hitting time of the extremum; Monte-Carlo simulation of these distributions using the characteristic functions[22, 21, 19, 20, 23]. Variations of the same technique for stable Lévy processes are used in [15, 17] to evaluate expectations of functions of $X_t$, $(X_t, \bar{X}_t)$, (joint) probability distributions in particular, and design Monte-Carlo simulation schemes for stable Lévy processes. The present paper demonstrates that the properties of SINH-regular, sSL- and SL-processes allow one to derive from the integral representation for the distribution of $\bar{X}_t$ asymptotic formulas. The new ingredient developed in the paper are novel representations of the Wiener-Hopf factors in terms of $\psi$, zeros of $q + \psi(\xi)$ and supports of the absolute continuous components of $\mathcal{G}_\pm$. We hope that these representations are of independent interest. The asymptotic coefficients are expressed in terms of the same quantities. We believe that the properties of SINH-regular and (s)SL processes and technique developed in the paper can be used to calculate the leading term of asymptotics of the distribution of $X$ conditioned on the survival, and many other situations, such as the calculation of of the limit of of the joint hit distribution of the running maximum of a Lévy process calculated in [27] in the BM model, and the asymptotics of the survival probability in the Sparre Andersen model in which the capital reserves are invested in a risky asset [26]. The representation of the asymptotic coefficient in terms of the first instantaneous moment and Stieltjes-Lévy measures can be used to study the contribution of different components of the risky asset to the survival probability in the long run.

The derivation of the asymptotic formulas uses the basic complex-analytical tools, essentially, the Cauchy integral formula only and resulting contour deformation technique. These tools suffice to evaluate integrals in many situations. ("The shortest path between two truths in the real domain passes through the complex domain." - J. Hadamard). The technique may seem not quite appropriate for specialists in probability. However, the results can be interpreted in terms of basic ingredients defining the process: Stieltjes-Lévy measures and characteristic exponent, therefore, the pure probabilist may try to prove the resulting formulas using purely probabilistic tools, and use the technique of the paper to derive formulas in new situations, and then prove them.

## References

[1] S. Asmussen. *Ruin Probabilities*. Number 2 in Advanced Series on Statistical Science and Applied Probability. World Scientific, River Edge, NJ, 2000.

[2] S. Asmussen, F. Avram, and M.R. Pistorius. Russian and American put options under exponential phase-type Lévy models. *Stochastic Processes and their Applications*, 109(1):79–111, 2004.

[3] F. Aurzada and T. Simon. Persistence probabilities & exponents. Working paper, March 2012. https://arXiv.org/abs/1203.6554.

[4] O.E. Barndorff-Nielsen. Processes of Normal Inverse Gaussian Type. *Finance and Stochastics*, 2:41–68, 1998.

[5] O.E. Barndorff-Nielsen and S.Z. Levendorskiĭ. Feller Processes of Normal Inverse Gaussian type. *Quantitative Finance*, 1:318–331, 2001.

[6] M. Boyarchenko, M. de Innocentis, and S. Levendorskiĭ. Prices of barrier and first-touch digital options in Lévy-driven models, near barrier. *International Journal of Theoretical and Applied Finance*, 14(7):1045–1090, 2011. Available at SSRN: http://papers.ssrn.com/abstract=1514025.

[7] M. Boyarchenko and S. Levendorskiĭ. Prices and sensitivities of barrier and first-touch digital options in Lévy-driven models. *International Journal of Theoretical and Applied Finance*, 12(8):1125–1170, December 2009.

[8] M. Boyarchenko and S. Levendorskiĭ. Ghost Calibration and Pricing Barrier Options and Credit Default Swaps in spectrally one-sided Lévy models: The Parabolic Laplace Inversion Method. *Quantitative Finance*, 15(3):421–441, 2015. Available at SSRN: http://ssrn.com/abstract=2445318.

[9] S. Boyarchenko and S. Levendorskiĭ. Generalizations of the Black-Scholes equation for truncated Lévy processes. Working Paper, University of Pennsylvania, April 1999.

[10] S. Boyarchenko and S. Levendorskiĭ. Option pricing for truncated Lévy processes. *International Journal of Theoretical and Applied Finance*, 3(3):549–552, July 2000.

[11] S. Boyarchenko and S. Levendorskiĭ. Barrier options and touch-and-out options under regular Lévy processes of exponential type. *Annals of Applied Probability*, 12(4):1261–1298, 2002.

[12] S. Boyarchenko and S. Levendorskiĭ. *Non-Gaussian Merton-Black-Scholes Theory*, volume 9 of *Adv. Ser. Stat. Sci. Appl. Probab.* World Scientific Publishing Co., River Edge, NJ, 2002.

[13] S. Boyarchenko and S. Levendorskiĭ. Perpetual American options under Lévy processes. *SIAM Journal on Control and Optimization*, 40(6):1663–1696, 2002.

[14] S. Boyarchenko and S. Levendorskiĭ. Sinh-acceleration: Efficient evaluation of probability distributions, option pricing, and Monte-Carlo simulations. *International Journal of Theoretical and Applied Finance*, 22(3):1950–011, 2019. DOI: 10.1142/S0219024919500110. Available at SSRN: https://ssrn.com/abstract=3129881 or http://dx.doi.org/10.2139/ssrn.3129881.

[15] S. Boyarchenko and S. Levendorskiĭ. Conformal accelerations method and efficient evaluation of stable distributions. *Acta Applicandae Mathematicae*, 169:711–765, 2020. Available at SSRN: https://ssrn.com/abstract=3206696 or http://dx.doi.org/10.2139/ssrn.3206696.

[16] S. Boyarchenko and S. Levendorskiĭ. Static and semi-static hedging as contrarian or conformist bets. *Mathematical Finance*, 3(30):921–960, 2020. Available at SSRN: https://ssrn.com/abstract=3329694 or http://arXiv.org/abs/1902.02854.

[17] S. Boyarchenko and S. Levendorskiĭ. Efficient evaluation of expectations of functions of a stable Lévy process and its extremum. Working paper, September 2022. Available at SSRN: http://ssrn.com/abstract=4229032 or http://arxiv.org/abs/2209.12349.

[18] S. Boyarchenko and S. Levendorskiĭ. Lévy models amenable to efficient calculations. Working paper, June 2022. Available at SSRN: https://ssrn.com/abstract=4116959 or http://arXiv.org/abs/2207.02359.

[19] S. Boyarchenko and S. Levendorskiĭ. Efficient evaluation of joint pdf of a Lévy process, its extremum, and hitting time of the extremum. Working paper, December 2023. Available at SSRN: http://ssrn.com/abstract=4656492 or http://arxiv.org/abs/2312.05222.

[20] S. Boyarchenko and S. Levendorskiĭ. Simulation of a Lévy process, its extremum, and hitting time of the extremum via characteristic functions. Working paper, December 2023. Available at SSRN: http://ssrn.com/abstract=4656466 or http://arxiv.org/abs/2312.03929.

[21] S. Boyarchenko and S. Levendorskiĭ. Efficient evaluation of double barrier options. *International Journal of Theoretical and Applied Finance*, 27(2):2450007, 2024. Available at SSRN: http://ssrn.com/abstract=4262396 or http://arxiv.org/abs/2211.07765.

[22] S. Boyarchenko and S. Levendorskiĭ. Efficient evaluation of expectations of functions of a Lévy process and its extremum. *Finance and Stochastics*, 2024. To appear. Working paper version is available at SSRN: https://ssrn.com/abstract=4140462 or http://arXiv.org/abs/2209.02793.

[23] S. Boyarchenko and S. Levendorskiĭ. Faster than fft. Proceedings icfda2024.sciencesconf.org 12th Conference on Fractional Differentiation and its Applications 9-12 July 2024, Bordeaux, France, October 2024.

[24] N. Cai and S.G. Kou. Option pricing under a mixed-exponential jump diffusion model. *Operations Research*, 60(1):64–77, 2012.

[25] P. Carr, H. Geman, D.B. Madan, and M. Yor. The fine structure of asset returns: an empirical investigation. *Journal of Business*, 75:305–332, 2002.

[26] E. Eberlein, Y. Kabanov, and T. Sshmidt. Ruin probabilities for a sparre andersen model with investment. *Stochastic Processes and their Applications*, 144:72–84, 2022.

[27] J.R. Furling, P. Salminen, and P. Vallois. On joint hit distribution of the running maximum of Brownian motion. *Stochastic Processes and Their Applications*, 150:1204–1221, 2022.

[28] S.G. Kou. A jump-diffusion model for option pricing. *Management Science*, 48(8):1086–1101, August 2002.

[29] A. Kuznetsov. Wiener-Hopf factorization and distribution of extrema for a family of Lévy processes. *Ann.Appl.Prob.*, 20(5):1801–1830, 2010.

[30] A. Kuznetsov, A.E. Kyprianou, and J.C. Pardo. Meromorphic Lévy processes and their fluctuation identities. *Annals of Applied Probability*, 22(3):1101–1135, 2012.

[31] S. Levendorskiĭ. *Degenerate Elliptic Equations*, volume 258 of *Mathematics and its Applications*. Kluwer Academic Publishers Group, Dordrecht, 1993.

[32] S. Levendorskiĭ. Pricing of the American put under Lévy processes. Research Report MaPhySto, Aarhus, 2002. Available at http://www.maphysto.dk/publications/MPS-RR/2002/44.pdf, http://www.maphysto.dk/cgi-bin/gp.cgi?publ=441.

[33] S. Levendorskiĭ. Pricing of the American put under Lévy processes. *International Journal of Theoretical and Applied Finance*, 7(3):303–335, May 2004.

[34] S. Levendorskiĭ. Convergence of Carr's Randomization Approximation Near Barrier. *SIAM FM*, 2(1):79–111, 2011.

[35] S. Levendorskiĭ. Method of paired contours and pricing barrier options and CDS of long maturities. *International Journal of Theoretical and Applied Finance*, 17(5):1–58, 2014. 1450033 (58 pages).

[36] S. Levendorskiĭ and B. Paneyakh. Degenerate elliptic equations and boundary value problems. In M.S. Agranovich and M.A. Shubin, editors, *Encyclopeadia of Mathematical Sciences, Vol.63*, pages 131–202. Springer, Berlin, 1994.

[37] A. Lipton. Assets with jumps. *Risk*, pages 149–153, September 2002.

[38] D.B. Madan and F. Milne. Option pricing with V.G. martingale components. *Mathematical Finance*, 1(4):39–55, 1991.

[39] D.B. Madan and M. Yor. Representing the CGMY and Meixner Lévy processes as time changed Brownian motions. *Journal of Computational Finance*, 12(1):27–47, 2009.

[40] R.C. Merton. Option pricing when underlying stock returns are discontinuous. *Journal of Financial Economics*, 3:125–144, 1976.

[41] J. Pitman and M. Yor. Infinitely divisible laws associated with hyperbolic functions. *Canadian Journal of Mathematics*, 55(2):292–330, 2003.

[42] K. Sato. *Lévy processes and infinitely divisible distributions*, volume 68 of *Cambridge Stud. Adv. Math.* Cambridge University Press, Cambridge, 1999.

[43] R. Schilling, R. Song, and Z. Vodrachěk. *Bernstein Functions: Theory and Applications*. De Gruyter, Berlin, 2d edition, 2012.

[44] W. Schoutens and J.L. Teugels. Lévy processes, polynomials and martingaless. *Communications in Statistics: Stochastic Models*, 14:335–349, 1998.